\documentclass[11pt]{amsart}
\usepackage{amssymb,amsthm}

\everymath{\displaystyle}

\theoremstyle{plain}
\newtheorem{theo}{Theorem}[section]
\newtheorem{prop}[theo]{Proposition}
\newtheorem{coro}[theo]{Corollary}
\newtheorem{lemma}[theo]{Lemma}
\newtheorem{teo}[theo]{Theorem}
\newtheorem{lem}[theo]{Lemma}
\newtheorem{cor}[theo]{Corollary}
\newtheorem*{fact}{Fact}
\theoremstyle{definition}
\newtheorem{defi}[theo]{Definition}

\newtheorem*{question}{Question}
\newtheorem*{claim}{Claim}

\theoremstyle{remark}
\newtheorem{rema}[theo]{Remark}

\newtheorem{remark}[theo]{Remark}

\newcommand{\fin}{\hfill$\square$}

\newcommand{\BB}{\mathbb{B}}

\newcommand{\DD}{\mathbb{D}}

\newcommand{\HH}{\mathbb{H}}

\newcommand{\NN}{\mathbb{N}}

\newcommand{\PP}{\mathbb{P}}

\newcommand{\RR}{\mathbb{R}}
\renewcommand{\SS}{\mathbb{S}}

\newcommand{\ZZ}{\mathbb{Z}}

\def\cB{{\mathcal B}}  \def\cH{{\mathcal H}}

\def\cU{{\mathcal U}} 
\def\cD{{\mathcal D}}   
 
\def\cE{{\mathcal E}}   
 
\def\cF{{\mathcal F}}

\newcommand{\tr}{\text{\rm tr}}
 
\newcommand{\II}{\text{\rm II}}
\newcommand{\la}{\langle}  \newcommand{\ra}{\rangle}
\newcommand{\AdS}{\mbox{AdS}}
\newcommand{\ADS}{\mathbb{ADS}}
\newcommand{\dS}{\mbox{dS}}
\newcommand{\DS}{\mathbb{DS}}
\newcommand{\wt}{\widetilde}

\newcommand{\Cl}{\mbox{Cl}}

\title[Cosmological time versus CMC time II]{Cosmological time versus
  CMC time II: the de Sitter and anti-de Sitter cases}

\author[L. Andersson]{Lars Andersson$^\star$}
\thanks{${}^\star$ Supported in part by the NSF, under contract no. DMS
0104402 with the University of Miami.}  \address{${}^\star$Albert
Einstein Institute, 
Am M\"uhlenberg 1, D-14476 Potsdam, Germany \and Department of Mathematics,
Univ. of Miami, Coral Gables, FL
 33124, USA}  
\email{larsa@math.miami.edu}

\author[T. Barbot]{Thierry Barbot$^\dagger$}

\thanks{${}^\dagger\; {}^{\ddagger}   \; {}^\S$ 
Supported in
part by ACI ``Structures g\'eom\'etriques et Trous Noirs''.}

\address{${}^\dagger$CNRS, UMPA, \'Ecole Normale Sup\'erieure de Lyon.}
\email{Thierry.BARBOT@umpa.ens-lyon.fr} 

\author[F. B\'eguin]{Fran\c cois B\'eguin$^\ddagger$}
\address{${}^\ddagger$Laboratoire de
Math\'ematiques, Univ. Paris Sud.}
\email{Francois.Beguin@math.u-psud.fr} 

\author[A. Zeghib]{Abdelghani Zeghib$^\S$}
 \address{${}^\S$CNRS, UMPA, \'Ecole
Normale Sup\'erieure de Lyon.}
\email{zeghib@umpa.ens-lyon.fr}

\keywords{}
\subjclass{}
\date{\today}

\begin{document}

\begin{abstract}
This paper continues the investigation of constant mean curvature (CMC)
time functions in maximal globally hyperbolic
spatially compact spacetimes of constant sectional curvature, which was
started in \cite{ABBZ}. In that paper, the case of flat spacetimes was
considered, and in the present paper, the remaining cases of 
negative curvature (i.e. anti-de Sitter) spacetimes and postitive curvature
(i.e. de Sitter) spacetimes is dealt with. 
As in the flat case, the existence of CMC time functions
is obtained by using the level sets of the cosmological time function as
barriers. A major part of the work consists of proving the required 
curvature estimates for these level sets.

The nonzero curvature case presents significant new difficulties, in part due
to the fact that the
topological structure of nonzero constant curvature spacetimes is much richer
than that of the flat spacetimes. Further, the timelike convergence
condition fails for de Sitter spacetimes, and hence uniqueness for 
CMC hypersurfaces fails in general. We characterize those
de Sitter spacetimes which admit CMC time functions (which are automatically unique), 
as well as those which admit CMC foliations but no CMC time function. 
\end{abstract} 

\maketitle


\section{Introduction} 

This paper is the second part of our investigation of constant mean
curvature time functions in maximal globally hyperbolic
spatially compact spacetimes of constant sectional   curvature.  
The first paper \cite{ABBZ} was devoted to the case
of flat spacetimes. The present paper concerns the remaining
cases,  namely spacetimes of positive constant curvature (de Sitter 
spacetimes), and  of negative constant curvature (anti-de Sitter 
spacetimes).  

The approach used in the present paper is
the same as that of \cite{ABBZ}. Thus, we shall study the properties of 
constant mean curvature 
time function by making use of
the cosmological time function, which is defined more directly in terms of
the spacetime geometry. To
achieve this, we need to understand the geometry of the levels of the
cosmological time functions of the spacetimes under consideration. Roughly
speaking, we must prove that each level of the cosmological time
function has almost constant mean curvature.  
 

The constant curvature spacetimes which shall be considered in this paper 
have locally trivial geometry, being
locally isometric to Minkowski space, de Sitter space or 
anti-de Sitter space and thus the partial differential equations aspect of the analysis of these
spacetimes is trivial. However, the topology of these spaces may be 
highly nontrivial and although the spacetimes under consideration have a local
isometry pseudo-group 
of maximal dimension, they typically have trivial (global)
isometry groups. 
Indeed, it is the interplay between the topology and the
causal structure of the spacetime which is the source of most of the
difficulties encountered in our work. 

In the flat case~\cite{ABBZ}, we could use known results on maximal globally hyperbolic 
flat spacetimes and their cosmological time functions,  due in particular to G. Mess (\cite{mess}) and  
F. Bonsante (\cite{bonsante}). Here, we have to prove similar results in the de Sitter
and anti-de Sitter cases.  This will lead us to prove many independent facts on the
geometry of \emph{domains of dependence} in anti-de Sitter space (\S\ref{sec:AdS-case} to \ref{proof.AdS}) and in de Sitter space (\S\ref{sec:dS-case} to \ref{proof.dS}).
Going from the flat case to the de Sitter and anti-de Sitter case is
not trivial.  Even if the local geometry is no less symmetric, the
global geometrical aspects are much harder to deal with.  The relation
between  the cases of de Sitter and anti-de Sitter spacetimes and the case of
flat spacetimes may be illustrated by considering  spherical and hyperbolic
geometry in comparison to  Euclidean geometry: the non-flat case presents many
additional difficulties.



Recall that a spacetime $(M,g)$ is said to be \emph{globally
  hyperbolic} if it admits a \emph{Cauchy hypersurface}, i.e. a
spacelike hypersurface $S$ which intersects every inextendible causal
curve at exactly one point. A globally hyperbolic spacetime is called
\emph{spatially compact} if its Cauchy hypersurfaces are compact. For
technical reasons, we will restrict ourselves to spatially
compact, maximal globally hyperbolic spacetimes (MGHC spacetimes for
short). Although this is a significant restriction, 
spacetimes of this class have been extensively studied,
especially as cosmological models. It is worth remarking that 
several authors, see eg. \cite{Bar1}, use the term 
cosmological spacetime to denote a MGHC spacetime satisfying the timelike
  convergence, or strong energy condition, 
i.e. $\text{Ric}(v,v)\geq 0$ for 
every timelike vector $v$. 
%
 
Among the spacetimes we consider here are those with positive constant
curvature, i.e. MGHC de Sitter spacetimes. The timelike convergence 
condition is violated in these spacetimes and hence the standard proof of
uniqueness of CMC foliations does not apply. Nevertheless, we shall
demonstrate the existence of a large class of MGHC de Sitter spacetimes which admit a
CMC time function, and thus a unique CMC foliation. 

%

The nonzero constant curvature spacetimes considered here are special cases
of spacetimes satisfying the vacuum Einstein equations with cosmological
constant. The current standard model of cosmology  has as an essential
element the accelerated expansion of the universe.  In order to achieve
accelerated expansion  the strong energy condition must be violated, which
leads one to consider spacetimes with positive cosmological constant,
i.e. spacetimes of de Sitter type. On the other hand, spacetimes of
anti-de Sitter type play an important role in the AdS/CFT correspondence,
which is currently being intensely investigated by string theorists.

\subsection{CMC time functions and CMC foliations.} \label{sec:CMCtime} 
We shall consider only time oriented spacetimes.  A globally hyperbolic
spacetime may be endowed with a time function,
i.e. a function $t:M\to \mathbb R$ which is strictly increasing on each
future directed causal curve. The trivial case is that of a direct metric
product $ M= (I, -dt^2)\times (N, h)$, where $I$ is an interval of ${\mathbb
R}$ and $(N, h)$ is a Riemannian manifold. In the general case,  a globally
hyperbolic spacetime  still has a topological product structure, but the
geometry may be highly distorted. It is attractive, from the mathematical as
well as the physical point of view, to  analyze the geometric distortion by
introducing a canonical time function, defined in a coordinate invariant
manner.  Introducing a canonical time function allows one  to describe the
spacetime as a one parameter family of Riemannian spaces indexed by time.
Here, we will consider CMC time functions (or CMC foliation when CMC time
functions do not exist).

In order to fix conventions, let  the second fundamental form of a spacelike
hypersurface $S$  be defined by $\II(X,Y) = \la \nu , \nabla_X Y\ra$ where
$\nu$ is the future oriented unit normal of $S$, and let the mean curvature
of $S$ be given by  $\tr(\II)/(n-1)$.

\begin{defi} 
Let $(M,g)$ be  time oriented spacetime.  A \emph{time function} on $M$ is a
function  $\tau: M \to \RR$ which is strictly increasing along any future
oriented causal curve.  A \emph{CMC time function} is a time function
$\tau_{cmc}:M \to \mathbb R$ such that the level $\tau_{cmc}^{-1}(a)$, if not
empty, is a Cauchy hypersurface with constant mean curvature $a$.
\end{defi}

\begin{defi}
A  \emph{CMC foliation} is a codimension one foliation whose leaves are
constant mean curvature spacelike hypersurfaces.
 \end{defi}
 
\begin{remark}
%
The existence of a  CMC time function is a considerably stronger condition
than the existence of a CMC foliation. In particular,  the definition of a
CMC time function requires not only  that the mean curvature of the
hypersurface $\tau_{cmc}^{-1}(a)$ is constant, but also that this mean
curvature is equal to $a$.  Hence, the mean curvature of the hypersurface
$\tau_{cmc}^{-1}(a)$ increases when $a$ increases. We do not require any
condition of this type for CMC foliations.

A consequence of the definition is that a CMC time function is \emph{always}
unique.  Actually, if a spacetime $M$ admits a CMC time function
$\tau_{cmc}$, then the foliation defined by the level sets of the function
$\tau_{cmc}$ is always the unique CMC foliation in $M$ (this is a
straightforward consequence of the maximum principle, see~\cite[\S
2]{BBZ}). Recall that, in general,  a  spacetime can admit infinitely many
CMC foliations.

Similarly, a CMC time function in a constant curvature MGHC spacetime is
automatically real analytic (see Proposition 5.12 of \cite{ABBZ}) whereas
this is not necessarily the case for CMC foliations (see
e.g. Proposition~\ref{p.foliations-de Sitter} and Remark~\ref{rk.feuillise}, item 2 and 3).
\end{remark}

As is well known, CMC hypersurfaces are solutions to a variational problem. 
There are deep connections between CMC hypersurfaces in both Riemannian
and Lorentzian spaces, and minimal surfaces, which are a classical subject
in differential geometry and geometric analysis. 
In general relativity, the CMC
time gauge plays an important role, and leads to a well posed Cauchy problem
for the Einstein equations. The CMC conjecture, one version of 
which may be formulated as
stating that a MGHC 
vacuum (i.e. Ricci flat) 
spacetime containing a CMC Cauchy hypersurface admits a global CMC time
function is one of the important
conjectures in general relativity, see \cite{andersson:survey} for
discussion. It should be noted, however, that there are spacetimes which
contain no CMC Cauchy surface. This was first pointed out by Bartnik
\cite{Bar1}. An example of a MGHC vacuum spacetime with this property was later
given by Chrusciel et al. \cite{chrusciel:etal:gluing}.

In spacetime dimension 3, the CMC time gauge leads naturally to a
formulation of the Einstein equations as a finite dimensional Hamiltonian
system on the cotangent bundle of Teichm\"uller space. See the introduction
to \cite{ABBZ} for further discussion.


\subsection{Statements of results} \label{sec:statements} 
%
%
Together with~\cite{ABBZ}, the present paper provides a complete
answer to the existence problem  of CMC time functions in the class of
MGHC spacetimes of {\it constant sectional curvature}.  

\subsubsection{The flat case}
We recall the main result of~\cite{ABBZ} (see also \cite{andflat,
  barflat}). 
\begin{theo}[\cite{ABBZ}]
\label{t.main-flat-case}
Let $(M,g)$ be a MGHC flat  spacetime. The following statements are true. 
\begin{enumerate}
\item If $(M,g)$ is not past (resp. future)
complete,
  then it admits a globally defined CMC time function $\tau_{cmc}:M\rightarrow
  I$ where $I=(-\infty,0)$ (resp. $I=(0,+\infty)$).    
\item If $(M,g)$ is causally 
complete then it admits a unique CMC foliation, but no globally defined CMC
time function. 
\end{enumerate}
\end{theo}

\subsubsection{The anti-de Sitter case}
The fact that the 
timelike convergence condition holds strictly in anti-de Sitter spacetimes (i.e. spacetimes with constant negative sectional curvature) simplifies the
analysis of CMC time functions. 
We shall prove the following result:

\begin{theo}[see  \S\ref{proof.AdS}]
\label{t.main-AdS-case}
Let $(M,g)$ be a MGHC spacetime with negative constant sectional
curvature. Then $(M,g)$ admits a globally defined CMC time function
$\tau_{cmc}:M\rightarrow (-\infty,\infty)$.
\end{theo}

\begin{rema}
Theorem~\ref{t.main-AdS-case}  was already proved in \cite{BBZ} in the
particular case where $\mbox{dim}(M)=3$.  The proof provided in~\cite{BBZ} uses
some sophisticated tools, such as the so-called \emph{Moncrief flow} on the
cotangent bundle of the Teichm\"uller space, which are very specific to the
case where $\mbox{dim}(M)=3$.
\end{rema}

\subsubsection{The de Sitter case}
In de Sitter spacetimes, i.e. spacetimes of constant  positive sectional
curvature, the timelike convergence condition fails to hold, and due to this
fact   the problem of existence of CMC time functions is most difficult in
this case.
%
Although they are quite delicate to deal with, MGHC de Sitter spacetimes are
very abundant and easy to construct. Any compact conformally flat Riemannian
manifold gives rise by means of a natural suspension process to a MGHC
de Sitter spacetime, and vice-versa. This classification is  essentially due
to K. Scannell (for more details, see section \ref{subsec.ds}). All of theses
spaces are (at least) future complete or past complete.

According to the nature of the  holonomy group of the associated conformally
flat Riemannian manifold, i.e.  the representation of its fundamental group
into  the M\"{o}bius group, MGHS de Sitter spacetimes split  into three types:
elliptic, parabolic and hyperbolic (reminiscent of the same classification in
Riemannian geometry).

Elliptic and parabolic de Sitter spacetimes admit a simple characterization.
\begin{itemize}
\item Every elliptic de Sitter spacetime is the quotient of the \emph{whole}
  de Sitter space by a finite group of isometries.
\item  Up to a finite cover, every parabolic dS spacetime is the quotient of
  some open domain of the de Sitter by a finite rank abelian group of
  isometries of parabolic type.
\end{itemize}
Using these geometrical descriptions, it is quite easy to prove that elliptic
and parabolic spacetimes do not admit any CMC time function, but admit CMC
foliations: More precisely, one has the following results:

\begin{prop}[see~\S\ref{sub.ellicmc}]
\label{t.elli-dS-case}
Let $(M,g)$ be an elliptic de Sitter MGHC spacetime. Then, $(M,g)$ admits no
CMC time function, but it admits (at least) a CMC foliation. More precisely:
\begin{enumerate}
\item if $(M,g)$ is isometric to the whole de Sitter space, it admits
  infinitely many CMC foliations.
\item if $(M,g)$ is isometric to a quotient of the de Sitter space by a
  non-trivial group, then there is a unique CMC foliation. Moreover, every
  CMC Cauchy hypersurface surface in $(M,g)$ is a leaf of this CMC foliation.
\end{enumerate}
\end{prop}

\begin{prop}[see \S\ref{sub.paracmc}]
\label{t.para-dS-case}
If $(M,g)$ is parabolic, then it  admits no CMC time function, but has a
unique CMC-foliation. Moreover, every CMC Cauchy surface in $(M,g)$ is
a leaf of this CMC foliation.
\end{prop}

``Most" de Sitter MGHC spacetimes are hyperbolic. Our last result, even if
non-optimal, tends to show that these spacetimes ``usually" admit CMC time
functions: 

\begin{theo}[see~\S\ref{sub.dscmc}]
\label{t.main-dS-case}
Let $(M,g)$ be a MGHC  hyperbolic de Sitter spacetime. After reversal of time, we
can assume that $M$ is future complete. Then, $(M,g)$ admits a partially
defined CMC time function $\tau_{cmc}:U\rightarrow I$ where $U$ is a
neighbourhood of the past end of $M$ and $I=(-\infty, \beta)$
for some $\beta \leq -1$. Moreover, $U$ is the whole spacetime $M$ and 
$\beta = -1$ in the following cases,
\begin{enumerate}
\item  $(M,g)$ has  dimension $2+1$,
\item  $(M,g)$ is a \textit{almost-fuchsian}, i.e. contains a Cauchy
  hypersurface with all principal curvatures  $< -1$.
\end{enumerate}
\end{theo}

\begin{remark}
Theorem~\ref{t.main-dS-case} is sharp in the following sense: for  any $n
\geq 4$, we will give examples of $n$-dimensional de Sitter MGHC spacetimes
which do not admit any global CMC time function  (see
section\ref{sss.no-CMC-time}).

A proof of Theorem~\ref{t.main-dS-case} in the particular case where
$\mbox{dim}(M)=3$ was given in \cite{Bar.Zeg}. This proofs relies on a
Theorem of F. Labourie on hyperbolic ends of $3$-dimensional manifolds, and
thus, is very specific to the $3$-dimensional case.
\end{remark}

\begin{remark}
There is a well-known natural duality between spacelike immersions
of hypersurfaces in de Sitter space and immersions of hypersurfaces in the hyperbolic space 
(see for example \cite[\S 5.2.3]{Bar.Zeg}).
This correspondance has the remarkable property to invert
principal curvatures: if $\lambda$ is a principal curvature of the spacelike
hypersurface immersed in de Sitter space, then the inverse $\lambda^{-1}$ is a principal curvature
of the corresponding hypersurface immersed in the hyperbolic space.

The notion of almost-fuchsian manifolds has been introduced by K. Krasnov
and J.-M. Schlenker in~\cite[\S 2.2]{KraSch} for the riemannian case. More precisely,
they defined \emph{almost-fuchsian hyperbolic\/} manifolds as
hyperbolic quasi-fuchsian manifolds containing a closed hypersurface $S$ with principal 
curvatures in $]-1, +1[$. 

For every $r>0$, let $S_r$ be the surface made of points at oriented distance $r$ from $S$.
Then, for $r$ converging to $-\infty$, the principal curvatures of $S_r$ all
tend to $-1$ (see~\cite[Lemma 2.7]{KraSch}). It follows that hyperbolic almost-fuchsian 
hyperbolic manifolds can be defined more precisely as hyperbolic quasi-fuchsian manifolds containing
a closed hypersurface $S$ with principal curvatures in $]-1, 0[$. 

Here we extended the notion of almost-fuchsian manifolds to the de Sitter case,
defining (future complete) almost-fuchsian de Sitter spacetimes
as MGHC de Sitter spacetimes containing a Cauchy hypersurface admitting principal curvatures
in $]-\infty, -1[$. It follows from the discussion above that this
terminology is consistent with respect to the Krasnov-Schlenker terminology
and the duality between de Sitter space and hyperbolic space.

Typical examples are \emph{fuchsian} spacetimes and small
deformations thereof 
(see Remark~\ref{rk.fuchsianalmost}).
\end{remark}

\section{Some general facts} 
\subsection{Cosmological time functions}
\label{s.cosmological-time}

In any spacetime $(M,g)$, one can define the \textit{cosmological time
  function\/}, see \cite{cosmic}, as follows:

\begin{defi}
The cosmological time function of a spacetime $(M,g)$ is the function
$\tau:M\rightarrow [0,+\infty]$ defined by
$$\tau(x)=\mbox{Sup}\{ L(c) \mid c \in {\mathcal R}^-(x) \},$$  where
${\mathcal R}^-(x)$ is the set of past-oriented  causal curves starting at
$x$, and $L(c)$ is  the lorentzian length of the causal curve $c$.
\end{defi}

This function  is in general badly behaved. For example,  in the case of
Minkowski space, the cosmological time function is everywhere infinite.

\begin{defi}
\label{d.regular}
A spacetime $(M,g)$ has \textit{regular cosmological time function} $\tau$ if
\begin{enumerate}
\item $M$ has \textit{finite existence time,\/} i.e. $\tau(x) < \infty$ for
   every $x$ in $M$,
\item for every past-oriented inextendible causal curve $c: [0, +\infty)
  \rightarrow M$, $\lim_{t \to \infty} \tau(c(t)) = 0$.
\end{enumerate}
\end{defi}

In \cite{cosmic}, Andersson, Galloway and Howard have proved that spacetimes
whose cosmological time function is regular enjoy many nice properties.

\begin{teo}
\label{teo.cosmogood}
If a spacetime $(M,g)$ has regular cosmological time function $\tau$, then
\begin{enumerate}
\item $M$ is globally hyperbolic,
\item $\tau$ is a time function, i.e. $\tau$ is continuous and is strictly
  increasing along future-oriented causal curves,
\item for each $x$ in $M$, there is a future-oriented timelike geodesic $c:
  (0, \tau(x)] \rightarrow M$ realizing the distance from the "initial
  singularity", that is, $c$ has unit speed, is maximal on each segment, and
  satisfies:
$$ c(\tau(x))) = x \mbox{ and } \tau(c(t)) = t \mbox{ for every }t
$$
\item $\tau$ is locally Lipschitz, and admits first and second derivative
almost everywhere.
\end{enumerate}
\end{teo}

\begin{remark}
\label{rk.reverse}
Similarly, for every spacetime $(M,g)$, one may define the \textit{reverse
cosmological time function} of $(M,g)$. This is the function $\widehat{\tau}:
M \rightarrow [0,+\infty]$ defined by
$$\widehat{\tau}(x)=\mbox{Sup}\{ L(c) / c \in {\mathcal R}^+(x) \},$$  where
${\mathcal R}^+(x)$ is the set of future-oriented  causal curves starting at
$x$, and $L(c)$ the lorentzian length of the causal curve $c$. Then one may
introduce the notion of spacetime with regular reverse cosmological time
function, and prove a result analogous to Theorem~\ref{teo.cosmogood}.
\end{remark}

\subsection{From barriers to CMC time functions}

In this section, for the reader convenience, we reproduce (more and less
classical) statements  on the  notions of \textit{generalized mean curvature}
and \textit{sequence of asymptotic barriers} as already presented in
\cite{ABBZ}.

For a $C^2$ strictly spacelike hypersurface $S$, let $\II$ and $H_S$ denote
the second fundamental form and mean curvature of $S$, respectively. These
objects were defined in section \ref{sec:CMCtime}.

\begin{defi}
\label{d.generalized-curvature}
Let $S$ be an edgeless  
achronal topological hypersurface in a spacetime $(M,g)$. We do not assume
$S$ to be differentiable.  Given a real number $c$, we will say that $S$
\textit{has generalized mean curvature bounded from above by $c$ at $x$},
denoted $H_S (x) \leq c$, if  there is a causally convex open neighborhood
$V$ of $x$ in $M$ and a smooth (i.e. $C^2$)   spacelike hypersurface
$\SS_x^-$ in $V$ such that
\begin{itemize}
\item $x\in\SS_x^-$ and $\SS_x^-$ is contained in the past of  $S \cap V$ (in
$V$),
\item the mean curvature of $\SS_x^-$ at $x$ is bounded from above by $c$.
\end{itemize} 

Similarly, we will say that $S$ \textit{has generalized mean curvature
bounded from below by $c$ at $x$}, denoted $H_S(x) \geq c$, if,   there is a
geodesically convex open neighborhood $V$ of $x$ in $M$ and a smooth
spacelike hypersurface $\SS_x^+$ in $V$ such that~:
\begin{itemize}
\item $x\in\SS_x^+$ and $\SS_x^+$ is contained in the past  of $S \cap V$
(with respect to   $V$),
\item the mean curvature of $\SS_x^+$ at $x$ is bounded from below by $c$.
\end{itemize}
We will write $H_S \geq c$ and $H_S \leq c$ to denote that $S$ has
generalized mean curvature bounded from below, respectively above, by $c$ for
all $x \in S$.
\end{defi}

\begin{defi}
Let $c$ be a real number. A \textit{pair of $c$-barriers} is a pair of
disjoint topological Cauchy hypersurfaces $(\Sigma^-,\Sigma^+)$ in $M$ such
that
\begin{itemize}
\item $\Sigma^+$ is in the future of $\Sigma^-$,
\item $H_{\Sigma^+} \leq c \leq H_{\Sigma^-}$ in the sense of definition
\ref{d.generalized-curvature}.
\end{itemize} 
\end{defi}

\begin{defi}
Let $\alpha$ be a real number.  A \textit{sequence of asymptotic past
$\alpha$-barriers} is a sequence of topological Cauchy hypersurfaces
$(\Sigma_m^-)_{m\in\NN}$ in $M$ such that
\begin{itemize}
\item 
$\Sigma_m^-$ tends to the past end of $M$ when $m\to +\infty$ (\textit{i.e.}
  given any compact subset $K$ of $M$, there exists $m_0$ such that $K$ is in
  the future of $\Sigma_m^-$ for every $m\geq m_0$),
\item 
$a_m^- \leq H_{\Sigma_m^-} \leq a_m^+$,  where $a_m^-$ and $a_m^+$ are real
numbers such that $\alpha<a_m^-\leq a_m^+$, and such that $a_m^+\rightarrow
\alpha$ when $m\to +\infty$.
\end{itemize}

Similarly, a \textit{sequence of asymptotic future $\beta$-barriers} is a
sequence of topological Cauchy hypersurfaces $(\Sigma_m^+)_{m\in\NN}$ in $M$
such that
\begin{itemize}
\item $\Sigma_m^+$ tends to the future end of $M$ when $m\to +\infty$,
\item $b_m^- \leq H_{\Sigma_m^+} \leq b_m^+$, where $b_m^-$ and $b_m^+$ are
  real numbers such that $b_m^-\leq b_m^+<b$, and such that $b_m^-\rightarrow
  \beta$ when $m\to +\infty$.
\end{itemize}
\end{defi}
Assume now that $(M,g)$ is an $n$-dimensional  MGHC   spacetime of  constant
curvature.

\begin{theo}
\label{t.foliation}
Assume that $(M,g)$  has constant curvature $k$, and admits a sequence of
asymptotic past $\alpha$-barriers and a sequence of asymptotic future
$\beta$-barriers. If $k\geq 0$, assume moreover that $(\alpha,\beta)\cap
[-\sqrt{k},\sqrt{k}]=\emptyset$. Then, $(M,g)$ admits a CMC-time
$\tau_{cmc}:M\rightarrow (\alpha,\beta)$.
\end{theo}

For the de Sitter case, we will also need the following intermediate (local)
statement (see Remark 5.11 in \cite{ABBZ}).

\begin{theo}
\label{t.half-theorem}
Assume that $(M,g)$ has constant curvature $k$ and  admits a sequence of
asymptotic past $\alpha$-barriers. If $k\geq 0$, assume moreover
$\alpha\notin [-\sqrt{k},\sqrt{k}]$. Then, $(M,g)$ admits a CMC time function
$\tau_{cmc}:U\rightarrow (\alpha,\beta)$ where $U$ is a neighbourhood of the
past end of $M$ (i.e. the past of a Cauchy hypersurface in $M$) and $\beta$
is a real number bigger than $\alpha$.  \fin
\end{theo}

\subsection{Spaces of constant curvature as $(G,X)$-structures}  
\label{s.(G,X)-structures}

Let $X$ be a manifold and $G$ be a group acting on $X$ with the following
property:  if an element $\gamma$ of  $G$ acts trivially on  an open subset
of $X$, then $\gamma$ is the identity element of $G$. A
\emph{$(G,X)$-structure} on  a   manifold   $M$ is  an   atlas
$(U_i,\phi_i)_{i\in I}$ where
\begin{itemize} 
\item $(U_i)_{i\in I}$ is a covering of $M$ by open subsets,
\item for every $i$, the map $\phi_i$ is a homeomorphism from $U_i$ to an
  open set in~$X$,
\item for every $i,j$, the transition map
$\phi_i\circ\phi_j^{-1}:\phi_j(U_i\cap U_j)\rightarrow\phi_i(U_i\cap U_j)$ is
the restriction of an element of $G$.
\end{itemize} 
Given a manifold $M$ equipped with a $(G,X)$-structure $(U_i,\phi_i)_{i\in
I}$, one can construct two important objects: a map $d:\widetilde
M\rightarrow  X$, called \emph{developing map}, and representation
$\rho:\pi_1(M)\rightarrow G$, called \emph{holonomy  representation}. The map
$d$ is a local homeomorphism (obtained by pasting together some lifts of the
$\phi_i$'s) and satisfies the following  equivariance property: for every
$\widetilde x\in\widetilde M$ and every $\gamma\in\pi_1(M)$, one has
$d(\gamma\cdot \widetilde x)=\rho(\gamma)\cdot d(\widetilde x)$. The map $d$
is unique up to post-composition by an element of $G$ (and the choice of $d$
obviously fully determines the representation $\rho$. In general, $d$ is
neither one-to-one, nor onto. A good reference for all these notions
is~\cite{Gol2}.

\bigskip

Now let $(M,g)$ be a $n$-dimensional spacetime with constant curvature $k=0$
(respectively $k=1$ and $k=-1$). Then it is well-known that every point in
$M$ admits a neighbourhood which is isometric to an open subset of the
Minkowski space $\mbox{Min}_n$ (respectively the de Sitter space $\mbox{dS}_n$
and  the anti-de Sitter space $\mbox{AdS}_n$). In other words, the lorentzian
metric on $M$ can be seen as a $(G,X)$-structure, where $X=\mbox{Min}_n$
(respectively $\mbox{dS}_n$ and $\mbox{AdS}_n$) and $G=\mbox{Isom}(X)$. Hence
the general theory provides us with a locally isometric developing map
$d:\widetilde M\to X$ and a representation $\rho:\pi_1(M)\to\mbox{Isom}(X)$
such that $d(\gamma\cdot \widetilde x)=\rho(\gamma)\cdot d(\widetilde x)$ for
every $\widetilde x\in\widetilde M$ and every $\gamma\in\pi_1(M)$. The map
$d$ and the representation $\rho$ will play a fundamental role in the proofs
of Theorems~\ref{t.main-AdS-case} and~\ref{t.main-dS-case}.

\section{Description of anti-de Sitter MGHC spacetimes} \label{sec:AdS-case} 
We now start our investigation of anti-de Sitter spacetimes.  Our goal is to
prove Theorem \ref{t.main-AdS-case}. According to Theorem~\ref{t.foliation},
this reduces to finding two sequences of asymptotic barriers. These sequences
of barriers will be provided by the levels of the cosmological time
function. Thus, we essentially need to prove curvature estimates for the
level sets of the cosmological time function of any anti-de Sitter MGHC
spacetime. A key point is that every MGHC spacetime with constant curvature
$-1$ is isometric to the quotient of a certain open domain in the
anti-de Sitter space $\mbox{AdS}_n$ by a discrete subgroup of
$\mbox{Isom}(\mbox{AdS}_n)$. A consequence is that studying the cosmological
time functions of anti-de Sitter MGHC spacetimes amounts to studying the
cosmological time functions of certain open domains in $\mbox{AdS}_n$. These
domains are called \emph{AdS regular domains}.

We will proceed as follows. In the present section,
we define \emph{AdS regular domains}, using the conformal structure of the
anti-de Sitter space. We shall also give two characterisation of AdS regular
domains, using the Klein model of the anti-de Sitter space. In
section~\ref{s.AdS-properties}, we  shall study the cosmological time and the
boundary of AdS regular domains. The desired estimates on the curvature of
the levels of the cosmological time of AdS regular domains will be obtained
in section~\ref{s.AdS-estimates}. Theorem~\ref{t.main-AdS-case} follows
easily from these estimates and from Theorem~\ref{t.foliation}.


\subsection{The linear model $\AdS_n$ of the anti-de Sitter space}
For $n\geq 2$, let $(x_1,\dots,x_{n+1})$ be the standard coordinates on
$\RR^{n+1}$, and consider the quadratic form
$Q_{2,n-1}=-x_1^2-x_2^2+x_3^2+\dots+x_{n+1}^2$. The \emph{linear model
$\AdS_n$ of the $n$-dimensional anti-de Sitter space} is the quadric
$(Q_{2,n-1}=-1)$, endowed with the lorentzian metric induced by $Q_{2,n-1}$.

It is very easy to see that $\AdS_n$ is diffeomorphic to
$\SS^1\times\DD^{n-1}$.  The geodesics of $\AdS_n$ are the connected
components of the intersections of $\AdS_n$ with the linear $2$-planes in
$\RR^{n+1}$. Similarly, the totally geodesic subspaces of dimension $k$ in
$\AdS_n$ are the connected components of the intersections of $AdS_n$ with
the linear subspaces of dimension $(k+1)$ in $\RR^{n+1}$.

A nice feature of the anti-de Sitter space is its simple conformal structure.

\begin{prop}
  \label{p.causal-structure}
The anti-de Sitter space $\AdS_n$ is conformally equivalent to
$(\SS^1\times\DD^{n-1},-dt^2+ds^2)$, where $dt^2$ is the standard riemannian
metric on $\SS^1=\RR/2\pi\ZZ$, where $ds^2$ is the standard metric (of
curvature $+1$) on the sphere $\SS^{n-1}$ and $\DD^{n-1}$ is the open upper
hemisphere of $\SS^{n-1}$.

Moreover, one can attach a Penrose boundary $\partial\wt\AdS_n$ to
$\wt\AdS_n$ such that $\wt\AdS_n\cup\partial\wt\AdS_n$ is conformally
equivalent to $(\SS^1\times\overline{\DD^{n-1}},-dt^2+ds^2)$, where
$\overline{\DD^{n-1}}$ is the closed upper hemisphere of $\SS^{n-1}$.
\end{prop}

Proposition~\ref{p.causal-structure} shows in particular that $\AdS_n$
contains many closed causal curves. One can overcome this difficulty by
considering the universal covering $\wt\AdS_n$ of $\AdS_n$. It follows from
Proposition~\ref{p.causal-structure} that $\wt\AdS_n$ is conformally
equivalent to $(\RR\times\DD^{n-1},-dt^2+ds^2)$, and admits a Penrose
boundary $\partial\wt\AdS_n$ such that  $\wt\AdS_n\cup\partial\wt\AdS_n$ is
conformally equivalent to   $(\RR\times\overline{\DD^{n-1}},-dt^2+ds^2)$. In
particular, $\wt\AdS_n$ and $\wt\AdS_n\cup\partial\wt\AdS_n$ are strongly
causal.

\begin{proof}[Proof of Proposition~\ref{p.causal-structure}]
See e.g. ~\cite[\S 4]{barbtz1} or~\cite[Proposition 4.16]{BBZ}.
\end{proof}

\subsection{AdS regular domains as subsets of $\AdS_n$}
In this paragraph, we will use the conformal completion
$\AdS_n\cup\partial\AdS_n$ of $\AdS_n$ to define the notion of \emph{AdS
regular domain}. Let us start by a remark.

\begin{rema}
\label{r.graph}
A subset $\wt\Lambda$ of $\partial\wt\AdS_n\approx
(\RR\times\SS^{n-2},-dt^2+ds^2)$ is achronal if and only if it is the graph
of a $1$-Lipschitz function $f: \Lambda_0 \rightarrow {\mathbb R}$ where
$\Lambda_0$ is a subset of ${\mathbb S}^{n-2}$ (endowed with its canonical
distance, induced by the metric $ds^2$ of curvature $1$).   In particular,
the achronal closed topological hypersurfaces in $\partial\wt\AdS_n$ are
exactly the graphs of the $1$-Lipschitz functions $f:\SS^{n-2}\to\RR$. In
particular, every closed achronal hypersurface in $\partial\wt\AdS_n$ is a
topological $(n-2)$-sphere.
\end{rema}

Let $\wt\Lambda$ be a closed achronal subset of $\partial\wt\AdS_n$, and
$\Lambda$ be the projection of $\wt\Lambda$ in $\partial\AdS_n$.  We denote
by  $\wt E(\wt\Lambda)$ the \emph{invisible domain} of $\wt\Lambda$ in
$\wt\AdS_n\cup\partial\wt\AdS_n$, that is,
$$\wt E(\wt\Lambda) = \left(\wt\AdS_n\cup\partial\wt\AdS_n\right) \setminus
\left(J^-(\wt\Lambda)\cup J^+(\wt\Lambda)\right)
$$ where $J^-(\wt\Lambda)$ and $J^+(\wt\Lambda)$ are the causal past and the
causal future of $\wt\Lambda$ in
$\wt\AdS_n\cup\partial\wt\AdS_n=(\RR\times\overline{\DD}^{n-1},-dt^2+ds^2)$.
We denote by $\Cl(\wt E(\wt\Lambda))$ the closure of $\wt E(\wt\Lambda)$ in
$\wt\AdS_n\cup\partial\wt\AdS_n$.  We denote by $E(\Lambda)$ the projection
of $\wt E(\wt \Lambda)$ in $\AdS_n\cup\partial\AdS_n$ (clearly, $E(\Lambda)$
only depends on $\Lambda$, not on $\wt\Lambda$).

\begin{defi}
A $n$-dimensional \emph{AdS regular domain} is a domain of the form
$E(\Lambda)$ where $\Lambda$ is the projection in $\partial\AdS_n$ of an
achronal topological $(n-2)$-sphere $\wt\Lambda\subset\partial\mbox{AdS}_n$.
\end{defi}

We will see later that regular domains satisfy several ``convexity
properties" (geodesic convexity, convexity in a projective space). The first
property of this kind concerns the causal structure.

\begin{rema}
For every closed achronal set $\wt\Lambda$ in $\partial\wt\AdS_n$, the
invisible domain  $\wt E(\wt \Lambda)$  is a \emph{causally  convex} subset
of $\wt\AdS_n\cup\partial\wt\AdS_n$: if $p,q\in \wt E(\wt\Lambda)$ then
$J^+(p)\cap J^-(q)\subset\wt E(\wt\Lambda)$, where $J^+(p)$ and $J^-(q)$ are
the causal past and future of $p$ and $q$ in
$\wt\AdS_n\cup\partial\wt\AdS_n$. This is an immediate consequence of the
definitions.
\end{rema}

The following remark is a key point for understanding the geometry of
AdS regular domains.

\begin{rema}
\label{r.f-f+}
Let $\wt\Lambda$ be a closed achronal subset of $\partial\wt\AdS_n$.
Recall that $\wt\Lambda$ is the graph of a $1$-Lipschitz function
$f:\Lambda_0\to\RR$ where $\Lambda_0$ is a closed subset of
$\SS^{n-2}$ (remark~\ref{r.graph}).  Define two functions
$f^-,f^+:\overline{\DD}^{n-1}\to\RR$ as follows: 
\begin{align*}
f^{-}(p) &= \mbox{Sup}_{q \in \Lambda_0} \{ f(q)-d(p,q) \} , \\
f^{+}(p) &= \mbox{Inf}_{q \in \Lambda_0} \{ f(q)+d(p,q) \} ,
\end{align*}
where $d$ is the distance induced by $ds^2$ on
$\overline{\DD}^{n-1}$. It is easy to check that 
$$
\wt E(\wt\Lambda)=\{(t,p)\in\RR\times\overline{\DD}^{n-1} \mid
f^-(p)<p<f^+(p)\}.
$$
\end{rema}

\begin{coro}
\label{c.one-to-one}
For every (non-empty) closed achronal set
$\wt\Lambda\subset\partial\wt\AdS_n$, the projection of $\wt
E(\wt\Lambda)$ on $E(\Lambda)$ is one-to-one. 
\end{coro}

\begin{proof}
We use the notations introduced in remark~\ref{r.f-f+}. For every
$p\in\overline{\DD}^{n-1}$, there exists 
a point $q\in\SS^{n-2}=\partial\overline{\DD}^{n-1}$ such that 
$d(p,q)\leq\pi/2$. Hence, for every $p\in\overline{\DD}^{n-1}$, we
have $f^+(p)-f^-(p)\leq\pi$. Hence $\wt
E(\wt\Lambda)$ is included in the set $
E=\{(t,p)\in\RR\times\overline{\DD^{n-1}} \mbox{ such that } f^-(p) <
t < f^-(p)+\pi\}. 
$
The projection of $\wt\AdS_n\cup\partial\wt\AdS_n=\RR\times\overline{\DD}^{n-1}$
on  $\AdS_n\cup\partial\AdS_n=(\RR/2\pi\ZZ)\times\overline{\DD}^{n-1}$ is
obviously one-to-one in restriction to $E$.
\end{proof}

\begin{coro}
\label{c.inter-boundary}
For every achronal topological $(n-2)$-sphere
$\wt\Lambda\subset\partial\wt\AdS_n$,
\begin{enumerate}
\item $\wt E(\wt\Lambda)$ is disjoint from $\partial\wt\AdS_n$
  (i.e. it is contained in $\wt\AdS_n$);
\item $\Cl\left(\wt E(\wt\Lambda)\right)\cap\partial\wt\AdS_n =
  \wt\Lambda$. 
\end{enumerate}
\end{coro}

\begin{proof}
We use the notations introduced in remark~\ref{r.f-f+}. 
Since $\wt\Lambda$ is a topological $(n-2)$-sphere, the set $\Lambda_0$
is the whole sphere $\SS^{n-2}$. Now observe that, for every
$p\in\SS^{n-2}=\Lambda_0$, one has $f^-(p)=f^+(p)=p$. Finally, recall
that $(t,p)\in\wt E(\wt\Lambda)$ (resp. $(t,p)\in\Cl(\wt
E(\wt\Lambda))$) if and only if $f^-(p)<t<f^+(p)$ (resp. $f^-(p)\leq
t\leq f^+(p)$). The corollary follows.
\end{proof}

The following notion will be useful later. 

\begin{defi}
Let $\Lambda_0$ be a closed subset of $\SS^{n-2}$, let
$f:\Lambda_0\to\RR$ be a 1-Lipschitz function, and 
$\wt\Lambda\subset\partial\wt\AdS_n$ be the graph of $f$. The achronal
set $\wt\Lambda$ is said to be \emph{pure lightlike} if  
$\Lambda_0$ contains two antipodal points $p_0$ and $-p_0$ 
on the sphere such that $f(p_0) = f(-p_0) +\pi$. 
\end{defi}

\begin{lem}
If $\wt\Lambda$ is pure lightlike, then $\wt E(\wt\Lambda)$ is empty.
\end{lem}

\begin{proof}
If $f:\Lambda_0\to\RR$ is $1$-Lipschitz, and if there exists two
antipodal points $p_0,-p_0\in\Lambda_0$ such that $f(p_0) = f(-p_0)
+\pi$, then it is easy to show that, for every element $p$ of
$\overline{\DD}^{n-1}$,  we have $f_-(p) = f_+(p) =
f(-p_0) + d(-p_0, p) = f(p_0) - d(p_0, p)$. The lemma follows.
\end{proof}

\subsection{The Klein model $\ADS_n$ of the anti-de Sitter
  space} 
We now consider the quotient $\SS(\RR^{n+1})$ of
$\RR^{n+1}\setminus\{0\}$ by positive homotheties. In other
words, $\SS(\RR^{n+1})$ is the double covering of the projective space
$\PP(\RR^{n+1})$. We denote by $\pi$ the projection of
$\RR^{n+1}$ on $\SS(\RR^{n+1})$. 
The projection $\pi$ is one-to-one in restriction to
$\AdS_n=(Q_{2,n-1}=-1)$. The \emph{Klein model} $\ADS_n$ of the
anti-de Sitter space is the projection of $\AdS_n$ in $\SS(\RR^{n+1})$,
endowed with the induced lorentzian metric. 

Observe that $\ADS_n$ is also the projection of the open domain of
$\RR^{n+1}$ defined by the inequality $(Q_{2,n-1}<0)$. It follows that
the topological boundary of $\ADS_n$ in $\SS(\RR^{n+1})$ is the
projection of the quadric $(Q_{2,n-1}=0)$; we will denote this
boundary by $\partial\ADS_n$. By construction, the projection $\pi$
defines an isometry between $\AdS_n$ and $\ADS_n$; one can easily
verify that this isometry can be continued to define a canonical
homeomorphism between 
$\AdS_n\cup\partial\AdS_n$ and $\ADS_n\cup\partial\ADS_n$. 

For every linear subspace $F$ of dimension $k+1$ in $\RR^{n+1}$, we
denote by $\SS(F)=\pi(F)$ the corresponding projective subspace of
dimension $k$ in $\SS(\RR^{n+1})$. The geodesics of $\ADS_n$ are the
connected components of the intersections of $\ADS_n$ with the
projective lines $\SS(F)$ of $\SS(\RR^{n+1})$. More generally, the
totally geodesic subspaces of dimension $k$ in $\ADS_n$ are the
connected components of the intersections of $\ADS_n$ with the
projective subspaces $\SS(F)$ of dimension $k$ of $\SS(\RR^{n+1})$. 

\begin{defi}
\label{def.affine}
An \emph{affine domain} of $\ADS_n$ is a connected component $U$ of
$\ADS_n\setminus\SS(F)$, where $\SS(F)$ is a projective hyperplane of 
$\SS(\RR^{n+1})$ such that $\SS(F)\cap\ADS_n$ is a spacelike
(totally geodesic) hypersurface. Let $V$ be the connected component of
$\SS(\RR^{n+1})\setminus\SS(F)$ containing $U$. The boundary $\partial
U\subset\partial\ADS_n$ of $U$ in $V$ is called the \textit{affine 
  boundary\/} of $U$. 
\end{defi}

\begin{remark}
\label{rk.affine}
Affine domains can be visualized in $\RR^{n}$. Indeed, let $U$ be an
affine domain in $\ADS_n$. By definition, there exists a a projective
hyperplane $\SS(F)$ in $\SS(\RR^{n+1})$ such that the hypersurface
$\SS(F)\cap\ADS_n$ is spacelike, and such that $U$ is one of the two
connected components of $\ADS_n\setminus\SS(F)$. We denote by $V$ the
connected component of $\SS(\RR^{n+1})$ containing $U$. Up to
composition by an element of the isometry group $SO_0(2,n-1)$ of
$Q_{2,n-1}$, we can assume that $\SS(F)$ is the projection of the
hyperplane $(x_1=0)$ in $\RR^{n+1}$ and $V$ is the projection of the
region $(x_1>0)$ in $\RR^{n+1}$. The map 
$$(x_1,x_2,\dots,x_{n+1})\mapsto
(u_1,\dots,u_n) :=
(\frac{x_2}{x_1},\frac{x_3}{x_1},\dots,\frac{x_{n+1}}{x_1})$$
induces a diffeomorphism between   
$V$ and $\RR^n$. In the coordinates
$(u_1,\dots,u_n)$, the image of  
the affine domain $U$ is to the region $(-u_1^2+u_2^2+\dots+u_n^2 <
1)$. The affine boundary $\partial U$ of $U$ corresponds to the
hyperboloid $(-u_1^2+u_2^2 =+\dots+u_n^2 = 1)$. The intersection of
$U$ with the totally geodesic subspaces of $\ADS_n$ correspond to the
intersections of the region $(-u_1^2+u_2^2+\dots+u_n^2 <
1)$ with the affine subspaces of $\RR^n$.
\end{remark}

\subsection{AdS regular domains as subsets of $\ADS_n$}
The canonical diffeomorphism between $\AdS_n\cup\partial\AdS_n$ and
$\ADS_n\cup\partial\ADS_n$ allows us to see AdS regular domains as
subsets of $\ADS_n$.
Nevertheless, it would be much more interesting to characterize AdS
regular domains directly as subsets of $\ADS_n$ without using the
identification of $\ADS_n\cup\partial\ADS_n$ with
$\AdS_n\cup\partial\AdS_n$; this is the purpose of the present
section. 
We start by stating the following lemma.

\begin{lem}
\label{l.affine-domain}
Let $\Lambda\subset\partial\AdS_n$ be the projection of a closed
achronal subset of $\partial\wt\AdS_n$ which is not pure lightlike. 
We see $\Lambda$ and $E(\Lambda)$ in $\ADS_n\cup\partial\ADS_n$. Then
$\Lambda$ and $E(\Lambda)$ are contained in the union $U\cup\partial
U$ of an affine domain and its affine boundary.
\end{lem}

\begin{proof} 
See \cite[Lemma~8.27]{barbtz1}. 
\end{proof}

Lemma~\ref{l.affine-domain} implies, in particular, that every AdS
regular domain is contained in an affine domain $U$ of $\ADS_n$. This
allows to visualize AdS regular domains as subsets of 
$\RR^n$ (see remark~\ref{rk.affine}). 

We will now use the pseudo-scalar product $\langle \cdot \mid \cdot \rangle$
associated with the quadratic form $Q_{2,n-1}$. It is important to
note that, although the real number $\langle x \mid y \rangle$ is
well-defined only for $x,y\in\RR^{n+1}$, the sign of $\langle x \mid y
\rangle$ is well-defined for $x,y\in\SS(\RR^{n+1})$. The following
lemma is easy but fundamental.

\begin{lemma}
\label{l.causally-related}
Let $U$ be an affine domain in $\ADS_n$ and $\partial
U\subset\partial\ADS_n$ be its affine boundary. Let $x$ be be a point
in $\partial U$, and $y$ be a point in $U\cup\partial U$. There exists
a causal (resp. timelike) curve joining $x$ to $y$ in $U\cup\partial
U$ if and only if $\langle x\mid y\rangle \geq 0$ (resp. $\langle
x\mid y\rangle>0$).  
\end{lemma}

\begin{proof}
See e.g. \cite[Proposition~5.10]{barbtz1} or~\cite[Proposition
4.19]{BBZ}. 
\end{proof}
 
Putting together the definition of the invisible domain $E(\Lambda)$
of a set $\Lambda\subset\partial\AdS_n$ and
Lemma~\ref{l.causally-related}, one easily proves the following.    

\begin{prop}
\label{pro.proj}
Let $\Lambda\subset\partial\AdS_n$ be the projection of a closed
achronal subset of $\partial\wt\AdS_n$ which is not pure lightlike.
If we see $\Lambda$ and $E(\Lambda)$ in the Klein
model $\ADS_n\cup\partial\ADS_n$, then
$$
E(\Lambda) = \{y \in \ADS_n\cup\partial\ADS_n\mbox{ such that }\langle
y  \mid x \rangle < 0 \mbox{ for every }x\in\Lambda\}).
$$
\end{prop}

\begin{remark}
\label{rk.nice}
A nice (and important) corollary of this Proposition is that the invisible
domain $E(\Lambda)$ associated with a set $\Lambda$ is
always geodesically convex, i.e. any geodesic joining two points in
$E(\Lambda)$ is contained in $E(\Lambda)$.  
\end{remark}

Proposition~\ref{pro.proj} provides a characterization of the AdS regular
domain associated to the projection of an achronal topological
$(n-2)$-sphere of $\partial\wt\AdS_n$. In order to obtain a complete
definition of AdS regular domains in $\ADS_n$, it remains to identify
the subsets of $\partial\ADS_n$ which corresponds to the projections
of achronal topological spheres $\partial\wt\AdS_n$. This is the
purpose of the following proposition, which easily follows from
Lemma~\ref{l.causally-related}.   

\begin{prop}
\label{le.list}
For $\Lambda\subset\partial\ADS_n$, the following assertions are
equivalent. 
\begin{enumerate}
\item when we see $\Lambda$ as a subset of $\partial\AdS_n$, it
  is the projection of an achronal subset of $\partial\wt\AdS_n$,
\item $\langle x \mid y \rangle$ is non-positive for every
  $x,y\in\Lambda$.
\end{enumerate}
Moreover, if $\Lambda$ satifies these assertions, $\Lambda$ is 
pure lightlike if and only it contains two antipodal points of
$\SS(\RR^{n+1})$. 
\end{prop}

Finally, we will give another characterization of the AdS regular
domains, using the duality for convex subsets of $\SS(\RR^n)$. 

Let us first recall some standard definitions. A \emph{convex cone}
$J$ of $\RR^{n+1}$ is a convex subset stable by positive homotheties. A
convex cone $J\subset \RR^{n+1}$ is said to be {\em proper} if it is
nonempty, and if its closure 
$\bar{J}$ does not contain a complete affine line. A \emph{convex
  subset} $C$ of $\SS(\RR^{n+1})$ is the projection of a convex cone
$J(C)$ of $\RR^n$ ; it is  \emph{proper} if $J(C)$ can be
chosen proper.  
Now, for any convex cone $J\subset\RR^{n+1}$, one can define the
\emph{dual} convex cone $J^*$ of $J$, 
$$
J^{\ast} = \{ x \in \RR^{n+1} \mbox{ such that } 
\langle x \mid y \rangle < 0 \mbox{ for all }y \in
\bar{J}\setminus\{0\}\}
$$
This allows one to associate a dual convex set $C^{\ast} \subset
\SS(\RR^{n+1})$ to any convex set $C\subset\SS(\RR^{n+1})$. Note
that $J^{\ast\ast} = J$ and $C^{\ast\ast} = C$. 

Using this duality, Proposition~\ref{pro.proj} can be reformulated as
follows.  
 
\begin{prop}
\label{lem.adsdual}
Let $\Lambda\subset\partial\AdS_n$ be the projection of a closed
achronal subset of $\partial\wt\AdS_n$. We see $\Lambda$ and
$E(\Lambda)$ in $\ADS_n\cup\partial\ADS_n$. Then the
domain $E(\Lambda)$ is the dual of the convex hull of $\Lambda$ in
$\SS(\RR^{n+1})$. \fin 
\end{prop}

In particular, AdS regular domains are the duals of the convex hulls
of the achronal topological $(n-2)$-sphere in $\partial\ADS_n$.

\subsection{Maximal globally hyperbolic spacetimes and regular
  domains.}  

The link between MGHC spacetimes with constant curvature $-1$ and
$\mbox{AdS}$ regular domains is made explicit by the following
theorem. 

\begin{theo}
\label{t.iso-quotient}
Every $n$-dimensional MGHC spacetime with constant curvature $-1$ is
isometric to the quotient of a regular domain in $\mbox{AdS}_n$ by a
torsion-free discrete subgroup of $\mbox{Isom}(\mbox{AdS}_n)$. 
\end{theo}

This result was proved by Mess in his celebrated preprint~\cite{mess}
(Mess only deals with the case where $n=3$, but his arguments also
apply in higher dimension). For the reader's convenience, we shall
recall the main steps of the proof (see~\cite[Corollary 
11.2]{barbtz1} for more details).   
  
\begin{proof}[Sketch of proof of Theorem~\ref{t.iso-quotient}]
Let $(M,g)$ be $n$-dimensional MGHC spacetime with constant curvature
$-1$. As explained in section~\ref{s.(G,X)-structures}, the theory of
$(G,X)$-structures provides us with a locally isometric developing
map $d:\widetilde M\to \mbox{AdS}_n$ and a holonomy representation
$\rho:\pi_1(M)\to\mbox{Isom}(\mbox{AdS}_n)$. Pick a Cauchy  
hypersurface $\Sigma$ in $M$, and a lift $\widetilde \Sigma$ of
$\Sigma$ in $\widetilde M$. Then $S:=d(\widetilde\Sigma)$ is an
immersed complete spacelike hypersurface in $\mbox{AdS}_n$. One can
prove that such a hypersurface is automatically properly embedded
and corresponds to the graph of a $1$-Lipschitz function
$f:\DD^2\to\SS^1$ in the conformal model
$(\SS_1\times\DD^2,-dt^2+ds^2)$. Such a function extends to a
$1$-Lipschitz function $\bar f$ defined on the closed disc
$\overline{\DD}^2$. This shows that the boundary $\partial S$ of $S$
in $\mbox{AdS}_n\cup\partial\mbox{AdS}_n$ is an achronal curve
contained in $\partial\mbox{AdS}_n$. 

On the one hand, it is easy to see that the Cauchy development
$D(S)$ coincides with the invisible domain $E(\partial S)$ (this
essentially relies on the fact that $S\cup\partial S$ is the graph of
a $1$-Lipschitz function, hence an achronal set in
$\widetilde{\mbox{AdS}}_n$).  
In particular, this shows that $D(S)$ is an AdS regular domain. 
 
On the other hand, one can prove that $M$ is isometric to the quotient
$\Gamma\setminus D(S)$, where $\Gamma:=\rho(\pi_1(M))$. Indeed, recall
that $S=d(\widetilde\Sigma)$  is a properly embedded
hypersurface. This shows that the group $\Gamma$ acts freely and
properly discontinuously on $S=d(\widetilde\Sigma)$. It is easy to
deduce that $\Gamma$ acts freely and properly discontinuously on the
Cauchy development $D(S)$. Hence the quotient $\Gamma\setminus D(S)$ is
a globally hyperbolic spacetime. Now, observe that $d(\widetilde M)$
is necessarly contained in $D(S)$ since $\widetilde\Sigma$ is a Cauchy
hypersurface in $\widetilde M$. Moreover, since $S$ is embedded in $M$, the
developing map $d$ is one-to-one in restriction  to
$\widetilde\Sigma$. It follows that $d$ is one-to-one on the Cauchy development of $\widetilde\Sigma$, i.e. on $\widetilde M$. Hence
the developing map $d$ induces an isometric embedding of $M$ in the
$\Gamma\setminus D(S)$. Since $M$ is maximal, this embedding must be
onto, and thus, $M$ is isometric to the quotient $\Gamma\setminus 
D(S)$.  
\end{proof}

\section{Cosmological time and horizons of AdS regular domains}
\label{s.AdS-properties} 

Throughout this section, we consider an achronal topological
$(n-2)$-sphere $\Lambda$ in $\partial\mbox{AdS}_n$, and the associated
AdS regular domain $E(\Lambda)$. 

\subsection{The cosmological time function} 

\begin{prop}
\label{pro.adsregular}
The AdS regular domain $E(\Lambda)$ has regular cosmological time.  
\end{prop}

\begin{proof} 
We recall that $\Lambda$ is, by definition, the projection of an
achronal topological sphere $\wt\Lambda\subset
\partial\wt\AdS_n$, and that $E(\Lambda)$ is the projection of
the invisible domain $\wt E(\wt\Lambda)$ of $\wt\Lambda$ in
$\wt\AdS_n\cup\partial\wt\AdS_n$. We will prove that $\wt
E(\wt\Lambda)$ has regular cosmological time. Since the projection of
$\wt E(\wt\Lambda)$ on $E(\Lambda)$ is one-to-one
(corollary~\ref{c.one-to-one}), this will imply that 
$E(\Lambda)$ also has regular cosmological time. We denote by $\wt\tau$ the
cosmological time of $\wt E(\wt\Lambda)$.

Let $x$ be a point in $\wt E(\wt\Lambda)$. On the one hand,
corollary~\ref{c.inter-boundary} states that $\Cl(\wt E(\wt\Lambda))$
is a compact subset of $\wt\AdS_n\cup\partial\wt\AdS_n$, and that $\Cl(\wt
E(\wt\Lambda))\cap \partial\wt\AdS_n=\wt\Lambda$. On the other hand,
since $x$ is in the invisible domain of $\wt\Lambda$, the set
$J^-(x)$ is disjoint from $\wt\Lambda$. Therefore $J^-(x)\cap\Cl(\wt    
E(\wt\Lambda))$ is a compact subset of $\wt\AdS_n$. Therefore
$J^-(x)\cap\Cl(\wt E(\wt\Lambda))$ is conformally equivalent to a
compact causally convex domain in $(\RR\times\DD^{n-1},-dt^2+ds^2)$ (with a
bounded conformal factor since everything is compact). It immediately
follows that the lengths of the past-directed causal curves starting
at $x$ contained in $\wt E(\wt\Lambda)$ is bounded (in other words,
$\wt\tau(x)$ is finite), and that, for every past-oriented
inextendible causal curve $c:[0,+\infty)\to \wt E(\wt\Lambda)$ with
$c(0)=x$, one has $\wt\tau(c(t))\to 0$ when $t\to\infty$.
This proves that $\wt E(\wt\Lambda)$ has regular cosmological time.
\end{proof}

Of course, since the definition of AdS regular domains is
``time-symmetric'', $E(\Lambda)$ also has regular reverse cosmological
time.

\subsection{Horizons}

According to Proposition~\ref{pro.adsregular} and
Theorem~\ref{teo.cosmogood}, $E(\Lambda)$ is globally 
hyperbolic. Hence its boundary in  
$\mbox{AdS}_n$ is a Cauchy horizon and enjoys all the  known
properties of Cauchy horizons (see for example \cite{beem}). In our
framework, this boundary is the union of two closed achronal subsets,
the past horizon ${\mathcal H}^-(\Lambda)$ and the future horizon 
${\mathcal H}^+(\Lambda)$. Observe that ${\mathcal H}^+(\Lambda)$
is in the future of ${\mathcal H}^-(\Lambda)$.

In the conformal model $(\DD^2\times\SS^1,-dt^2+ds^2)$, the horizons
${\mathcal H}^-(\Lambda)$ and ${\mathcal H}^+(\Lambda)$ are
the graphs of the functions $f^+$ and $f^-$ defined in
remark~\ref{r.f-f+}. In the Klein model, $E(\Lambda)$ is a convex domain,
and the union ${\mathcal H}^-(\Lambda)\cup {\mathcal H}^+(\Lambda)$ is
the topological boundary of this convex domain. We can therefore consider
support hyperplanes to $E(\Lambda)$ at some point $p\in {\mathcal
  H}^\pm(\Lambda)$. These are projective hyperplanes in
$\SS(\RR^{n+1})$. It is quite clear that, for such a support 
hyperplane $H\subset\SS(\RR^{n+1})$, the corresponding totally geodesic
hypersurface $H\cap\ADS_n$ is degenerate or spacelike (otherwise,
$H$ would intersect transversally the achronal hypersurface ${\mathcal
  H}^\pm(\Lambda)$, and this would contradict the fact 
that $H$ is a support hyperplane of $E(\Lambda$.

The following is the analogous of Lemma 3.1 in \cite{ABBZ}.

\begin{prop}
\label{lem.adshori}
Let $p$ a point of the past horizon ${\mathcal H}^-(\Lambda)$ of
$E(\Lambda)$. Let $C(p) \subset T_p\mbox{AdS}$ be the set of the
future directed unit tangent vectors orthogonal to the support hyperplanes
of $E(\Lambda)$ at $p$. Then: 
\begin{enumerate}
\item $C(p)$ is the convex hull of its lightlike elements. 
\item If $c$ is a future complete geodesic ray starting at $p$ whose
  tangent vector at $p$ is a lightlike element of $C(p)$, then the
  future endpoint of $c$ is in $\Lambda$. 
\end{enumerate}
\end{prop}

\begin{proof}
First of all, we need to understand better the link between the way
the elements of $C(p)$ are associated to the support planes of
$E(\Lambda)$ at $p$.
Let $H$ be a support hyperplane of $E(\Lambda)$  at $p$. Then
$H=\SS(u^\perp)$ where $u$ is an element of $\RR^{n+1}$ such that
\begin{enumerate}
\item[(i)] $\langle u \mid u\rangle\leq 0$ (since $H=\SS(u^\perp)$ is
  spacelike or lightlike); 
\item[(ii)] $\langle p \mid u \rangle = 0$ (since $p\in
  \SS(u^\perp)$);
\item[(iii)] $\langle x \mid u \rangle \leq 0$ for every $x\in
  E(\Lambda)$ (since $H=\SS(u^\perp)$ is a support hyperplane of
  $E(\Lambda)$, and since, up to replacing $u$ by $-u$, we can assume
  that $u$ and $E(\Lambda)$ are on the same side of $H$). 
\end{enumerate}
Observe that this property~(i) implies that the
  projection $[u]$ of $u$ in $\SS(\RR^n)$ belongs to
  $\ADS_n\cup\partial\ADS_n$. Also observe that $[u]$ and
  $E(\Lambda)$ being on the same side of $H$, the point $[u]$ must be
  in the future of $p$. Consider the 
$2$-plane $P_u$ containing $u$ and $p$. The projection of $\SS(P_u)$
of $P_u$ is a causal geodesic $\gamma_u$ containing $p$ and
orthogonal to $H$. If $[u]\in\ADS_n$, then $[u]\in\gamma_u$; if
$[u]\in\partial\ADS_n$, then $[u]$ is the final extremity of
$\gamma_u$. We will denote by $v_u$ be the future directed unit tangent
vector of $\gamma_u$ at $p$. 

The set $C(p)$ is the set of all the vectors $v_u$ when
$H=\SS(u^\perp)$ ranges other the set of all the support
hyperplanes of $E(\Lambda)$ at $p$. It is important to note 
that $v_u$ is lightlike if and only if $u$ is lightlike, i.e. 
if and only if $H=\SS(u^\perp)$ is a lightlike hyperplane.

Now we will prove item~(1). For this purpose, let us consider a
support plane $H=S(u^\perp)$ of $E(\Lambda)$ at $p$. We know that
$\langle x \mid u \rangle \leq 0$ for every  $x\in E(\Lambda)$. We
also know that 
$E(\Lambda)$ is the dual of the convex hull of  $\Lambda$ in
$\SS(\RR^n)$ (Proposition~\ref{lem.adsdual}). This implies that the projection
$[u]\in\SS(\RR^n)$ of $u\in\RR^n$ belongs to the convex hull in
$\SS(\RR^n)$ of $\Lambda$. Hence, we can write 
$u$ as a convex combination $u = \sum a_iu_i$ where the $u_i$'s are
elements of $\RR_n$ projecting onto elements of $\Lambda$ and the
$a_i$ are positive number (equivalently, $v_u$ is a convex combination
of the $v_{u_i}$'s). We know that the scalar $\langle p \mid u \rangle$
is equal to zero: $\sum a_i \langle p \mid u_i \rangle =0$. But 
all the terms of this sum are nonpositive. Therefore
$\langle p \mid u_i \rangle=0$ for every $i$. As a consequence,
$S(u_i^\perp)$ is a support plane of $E(\Lambda)$ at $p$
for every $i$ (equivalently, $v_{u_i}$ is an element of $C(p)$ for every
$i$). Moreover, $H_i=\SS(u_i^\perp)$ is a lightlike hyperplane for every $i$
(equivalently, $v_{u_i}$ is lightlike for every $i$). So, we have proved
that $v_u$ is a convex combination of elements lightlike elements
$C(p)$. This completes the proof of~(1).  

It remains to prove item~(2). For this purpose, we consider a support
plane $H=S(u^\perp)$ of $E(\Lambda)$ at $p$, and the associated element
$v_u$ of $C(p)$. We assume that $H$ is lightlike (equivalently that
$v_u$ is lightlike). 
Just as above, we write $u = \sum a_iu_i $ where the $u_i$'s
projecting on elements of 
$\Lambda$, and the $a_i$'s are positive. By hypothesis, the norm of
$u$ is equal to zero: $\sum a_ia_j \langle u_i \mid u_j
\rangle=0$. But, according to Proposition~\ref{le.list}, the scalar product
$\langle u_i \mid u_j \rangle$ is non-positive for every $i,j$.
Hence, $\langle u_i \mid u_j \rangle$ must be equal zero for every
$i,j$. Hence, the  subspace $F$ spanned by the $u_i$'s is (totally)
isotropic, which implies it  is either $1$-dimensional or $2$-dimensional. 
In the first case, $[u_i]=[u_j]$ for all $i$, $j$, and in the second one,
$S(F)$ is a lightlike geodesic containing all the $[u_i]$'s. In both
cases, we deduce that $[u]$ belongs to the segment joining $[u_i]$ to
$[u_j]$ for some $i,j$. It follows that $[u]$ belongs to $\Lambda$
(since $\Lambda$ is achronal, every lightlike segment with both ends
in $\Lambda$ is contained in $\Lambda$). Now recall that $v_u$ is the
tangent vector at $p$ of the geodesic segment joining $p$ and
$[u]$. Hence, we have proved that the future extremity of the
lightlike ray starting at $p$ with tangent vector $v_u$ is in
$\Lambda$. This completes the proof of~(2). 
\end{proof}

\begin{remark}
Of course, a similar statement holds for the future horizon 
${\mathcal H}^+(\Lambda)$ but where complete null rays 
contained in the horizon are now past oriented.
\end{remark}

\subsection{Retraction onto the horizon}

According to point $(3)$ in Theorem~\ref{teo.cosmogood}, for every point $x$ 
in the regular domain, there exists at least one maximal timelike
geodesic ray with future endpoint  $x$ realizing the ``distance to the initial 
singularity'': we call such a geodesic ray a \textit{realizing\/}
geodesic for $x$. 

\begin{defi}
The region $\{ \tau < \pi/2 \}$ of the AdS regular domain $E(\Lambda)$
is denoted $E_0^-(\Lambda)$ and called the {\bf past tight} region of
$E(\Lambda)$. 
\end{defi}

\begin{prop}
\label{pro.tout}
Let $x$ be an element of the past tight region $E_0^-(\Lambda)$ of
$E(\Lambda)$. Then, there is an unique realizing geodesic for $x$.     
\end{prop}

This proposition means that the past tight region is foliated by inextendible 
timelike geodesics on which $\tau$ restricts as a unit speed parameter.

\begin{proof}
Consider an affine domain $U$ containing $E(\Lambda)$ (see
Proposition~\ref{l.affine-domain}). 
In some coordinate system $(u_1, u_2, \ldots , u_n)$ the domain  
$U$ is the region $\{ -u_1^2 + u_2^2 + \ldots + u_n^2 < 1 \}$, and $x$ has 
zero coordinates (see definition~\ref{def.affine} and
remark~\ref{rk.affine}). 
Initial extremities of realizing geodesics for $x$ are points $z$ in 
${\mathcal H}^-(\Lambda)$ such that $d(x,z) = \tau(x)$, where $d(x,z)$
is the  
time length
of a past oriented timelike geodesic in 
${\mathbb A}{\mathbb D}{\mathbb S}_n$ starting from 
$x$ and
ending to $z$ (hence $d(x,z) = 0$ if $z$ is not in the past of
$x$). For each $\tau$,  we have: 
${\mathcal E}_\tau = \{ z \in U  | d(x,z) \geq \tau \} =  
\{ -u_1^2 + u_2^2 + \ldots + u_n^2 \leq -\tan^2(\tau), x_1 < 0  \}$. 
If $\tau < \tau'$, 
then ${\mathcal E}_{\tau'} \subset 
{\mathcal E}_\tau$. Since $E(\Lambda)$ is causally convex, one has:
$$
\tau(x) =  \sup\{ \tau  | {\mathcal E}_\tau \cap {\mathcal H}_-(\Lambda) \neq 
\emptyset \} 
$$
Let $y$, $y'$ be initial extremities of realizing geodesics for $x$:
they both belong  to ${\mathcal E}_{\tau(x)} \cap E(\Lambda)$. Assume
by contradiction that $y \neq y'$, and take any element $z$ in the
interior of the segment $[y,y']$. On the one hand, since $E(\Lambda)$
is geodesically convex, $z$ belongs to $E(\Lambda)$. 
On the other hand, $z$ belongs to the interior of ${\mathcal E}_{\tau(x)}$ 
(since the hyperboloid  
$\{ -u_1^2 + u_2^2 + \ldots + u_n^2 = -\tan^2(\tau(x)), x_1 < 0  \}$
is concave). 
Hence, the length of the geodesic segment $[x,z]$ is strictly bigger
than $\tau(x)$. Contradiction. 
\end{proof}


\begin{prop}
\label{pro.tight}
Let $c: (0, T] \rightarrow E_0^-(\Lambda)$ be a future
oriented timelike geodesic whose initial extremity  $p:=\lim_{t\to 0}c(t)$ is in the past
horizon $\cH^-(\Lambda)$. Then the following assertions are equivalent.
\begin{enumerate}
\item For every $t\in(0,T]$, $c_{|[0,t]}$ is a realising geodesic for the
  point $c(t)$.
\item There exists $t_0\in (0, T]$ such that $c((0,t_0])$ is a realizing
  geodesic for the point $c(t)$.
\item $c$ is orthogonal to a support hyperplane of $E(\Lambda)$ at
  $p:=\lim_{t\to 0}c(t)$.  
\end{enumerate}
\end{prop}

\begin{proof}
Obviously (1)$\Rightarrow$(2).

Assume that there exists $t_0\in (0, T]$ such that $c((0,t_0])$ is a
realizing geodesic for the point $c(t_0)$. Let $x:=c(t_0)$ and 
  $p:=\lim_{t\to 0}c(t)$. The level 
  set $\{ z  | d(x,z) =\tau(x)\}$ is a smooth hypersurface 
in ${\mathbb A}{\mathbb D}{\mathbb S}_n$
and its tangent space at $p$ is the orthogonal in 
$T_p{\mathbb A}{\mathbb D}{\mathbb S}_n$ of
the vector tangent to $c$. If this tangent space is not tangent to a 
support hyperplane of 
${\mathcal H}^-(\Lambda)$ then, the set $\{ d(x, z) =\tau(x) \}$  intersects 
$E(\Lambda)$.  This  would imply
that $p$ is not a minimum point  for the restriction
of $d(x, .)$ to $E(\Lambda)$. This is a contradiction since the
  restriction of $c((0,t_0])$ is a realizing geodesic for $x$. Hence
  (2)$\Rightarrow$(3). 

Now assume that $c$ is orthogonal to a support hyperplane 
of ${\mathcal H}^-(\Lambda)$ at $p$, and let $x$ be a point of  $c$.
Consider an affine domain $U$ centered at $x$.
The hyperboloid $\{ z \in U  | d(x,z) = d(x,p) \}$ is orthogonal to $c$ at $p$:
hence, by hypothesis, its tangent space at $p$ is a support hyperplane of  
$E(\Lambda)$.
Since ${\mathcal H}^-(\Lambda)$ is convex whereas the hyperboloid is
strictly concave, the intersection of
${\mathcal E}(\Lambda)$
with  ${\mathcal H}^-(\Lambda)$ is $\{ p\}$. This means that $p$ is a
minimum point for  $d(x,.)$. Therefore, $[x,y]$ is a realizing
geodesic for $x$. Hence (3)$\Rightarrow$(1). 
\end{proof}

\begin{rema}
Using the reverse cosmological time $\widehat\tau$ instead of $\tau$,
one can define the \emph{future tight region} $E_0^+(\Lambda)$ of
$E(\Lambda)$, and prove some analogs of Propositions~\ref{pro.tout}
and~\ref{pro.tight}. 
\end{rema}

\section{AdS regular domains: curvature estimates of cosmological levels}
\label{s.AdS-estimates}

We are now able to state the main result on curvature estimates of
the level sets of the cosmological time of an AdS regular domain.

\begin{teo}
\label{t.barriers-AdS}
Let $E^-(\Lambda)$ be the past tight region of an AdS regular domain, and
$\tau: E_0^-(\Lambda) \rightarrow (0, \pi/2)$  be the associated cosmological
time. For every $a \in (0, \pi/2)$, the generalized mean curvature of
the level set $S_a = \tau^{-1}(a)$ satisfies 
$$ -\cot(a) \leq H_{S_{a}} \leq -\frac{1}{n-1}\cot(a) +
  \frac{n-2}{n-1}\tan(a).$$ 
%
%
\end{teo}

\begin{proof}
Let $x$ be a point on the level set $S_a$. We denote by $c: [0,
a] \rightarrow E^-(\Lambda)$ the unique realizing geodesic for $x$,
with initial extremity $p = r(x)$. Let $v$ be the future oriented unit
speed tangent vector of $c$ at $p$. We denote as before $C(p)$ the set of
vectors in $T_p\ADS_n$ orthogonal to support hyperplanes of the past
horizon at $p$.   Our goal is to construct two local surfaces
${\mathbb S}_x^+$,  ${\mathbb S}_x^-$ containing $x$, respectively in
the future and the past of $S_{a}$, and with known mean curvature at $x$
(recall Definition~\ref{d.generalized-curvature}).

\medskip

\noindent {\it Construction of ${\mathbb S}_x^+$.} 
The construction of the upper barrier ${\mathbb S}^+_x$ is similar to
the construction in the flat case: take a portion near $x$ of the set
of points at lorentzian  distance $a$ from $p = r(x)$. The mean
curvature of ${\mathbb S}^+_x$ is $-\cot(a)$, its tangent hyperplane
at $x$ is the hyperplane orthogonal to $c$ at $x$. 

\medskip

\noindent {\it Construction of ${\mathbb S}_x^-$.}
Let $\widetilde\Lambda$ be a lift of $\Lambda$ in
$\partial\wt\AdS_n\simeq \RR\times\SS^{n-2}$. We recall
that $\widetilde\Lambda$ can be seen as the graph 
 of a $1$-Lipschitz function $f:\SS^{n-2}\to\RR$.
By Proposition~\ref{pro.tight}, the vector $v$ is in $C(p)$. 
Hence, Proposition~\ref{lem.adshori} implies that
there is a finite set $\{v_1,\dots,v_l\}$ of lightlike elements of
$C(p)$ such that $v$ is in the convex hull of $\{v_1,\dots,v_l\}$.
According to Proposition~\ref{lem.adshori} the future extremities 
of the lightlike geodesics whose tangent vectors at $p$ are
$v_1,\dots,v_l$ belong to $\Lambda$. Let $B$ be the finite subset of
$\Lambda$ made of these future extremities, and $\wt B$ the
corresponding subset of $\wt\Lambda$.
Then $\wt B$ is the graph of a $1$-Lipschitz  function $f_B: B_0
\rightarrow \mathbb R$ where $B_0$ 
is a finite subset of ${\mathbb S}^{n-2}$ (see remark~\ref{r.graph}). Let
$\wt\Lambda_B$ be the graph of the $f_B^-: {\mathbb S}^{n-2} \rightarrow
\mathbb R$ defined remark~\ref{r.f-f+}, and $\Lambda_B$ be the projection of
$\wt\Lambda_B$. We define our  hypersurface ${\mathbb S}_{x}^{-}$ to be the $a$-level set of the cosmological time of the domain $E(\Lambda_B)$. 

Let us check that $\SS_x^-$ satisfies the required properties: $x\in\SS_x^-$ and $\SS_x^-$ is in the past of $S_a$. Since $\Lambda_B$ subset $\Lambda$, the invisible domain $E(\Lambda_B)$ contains the invsible domain $E(\Lambda)$, and hence the hypersurface $\SS_x^-$ is in the past of the hypersurface $S_a$. For each $x\in\Lambda_B$, there is a future directed lightlike geodesic ray starting at $p$ whose endpoint is equal to $x$. It follows that $p\in\cH^-(\Lambda_B)$. By contruction, the vectors $v_1,\dots,v_l$ are orthogonal to support hyperplanes of $E(\Lambda_B)$ at $p$. Hence $v\in\mbox{Conv}(v_1,\dots,v_l)$ is also  orthogonal to a support hyperplane of $E(\Lambda_B)$ at $p$. According to Proposition~\ref{pro.tight}, this implies that $c$ is a realizing geodesic  in $E(\Lambda_B)$. It follows that $x=c(a)$ belongs to the $a$-level set of the cosmological time of $E(\Lambda_B)$, i.e. $x\in\SS_x^-$.

We are left to evaluate the mean curvature of the hypersurface $\SS_x^-$ at $x$.
The finite set $B$ is the projection of a set $\widehat{B}$
of null vectors in $E_n$. Let $F$ be the vector space spanned by $\widehat{B}$,
and let $F^\perp$ be the subspace orthogonal to $F$. Let $1+d$ be the dimension
of $F$. The convex hull of $\widehat{B}$ contains a timelike element $\hat{q}$
with $Q_{2,n-1}$-norm $-1$: the dual to the spacelike support hyperplane at 
$p$ orthogonal to $v$. This point $\hat{q}$ can also
be defined as the unique element of $\mbox{AdS}_n$ 
projecting on $q = c(\pi/2)$.

Similarly, $F^\perp$ contains a timelike vector:
the lift $\hat{p}$ in $E_n$ of $p$, let us  say,
$Q_{2,n-1}(\hat{p}) = -1$.
It follows that $F \cap F^\perp = \{ 0 \}$, $F$ has signature $(1,d)$,
and $F^\perp$ has signature $(1,n-d-1)$. 

Let $G \approx \mbox{SO}_0(1,
n-d-1)$ be the subgroup of $\mbox{SO}_0(2, n-1)$ made of the elements acting
trivially on $F$. The group $G$ preserves $\widehat{B}$.  
It follows that its induced action on $S(E)$ preserves
$E(\Lambda_B)$. This action preserves the cosmological time $\tau_B$
of $E(\Lambda_B)$. The $G$-orbit of $p$ is a connected component of 
the  geodesic subspace $S(F^\perp) \cap {\mathbb A}{\mathbb
  D}{\mathbb S}_n$. 

Let $F_1$ be the subspace $F^\perp \oplus \langle
\hat{q} \rangle$.  Observe that $\hat{q}$ is a fixed point for the
action of $G$. 
The projection  $A_1$ of $F_1 \cap \mbox{AdS}_n$ in $S(E_n)$ is a copy
of the Klein model of the anti de Sitter space of dimension $n-d$. It
contains $x$ which is the projection of $\hat{x} = \cos(a)\hat{p}+
\sin(a)\hat{q}$.  The $G$-orbit of $x$ is contained in the
cosmological level $\tau_B^{-1}(a)$. On the other hand, this $G$-orbit
in the anti de Sitter space $A_1$ is the set of initial extremities of
future oriented timelike geodesics with future extremity $q$ and of
length $\pi/2 - a$. Hence, it is an umbilical submanifold with principal
curvatures $\cot(\pi/2 - a) = \tan(a)$. This $G$-orbit is orthogonal
to $r^{-1}(p)$, and in $r^{-1}(p) \subset S(F)$,  the cosmological 
time $\tau_B$ is simply the lorentzian distance to $p$:
$\tau_B^{-1}(a) \cap r^{-1}(p)$ is an umbilical submanifold with
principal curvatures $-\cot(a)$. Hence, the mean curvature of
${\mathbb S}_x^-=\tau_B^{-1}(a)$ at points in $r^{-1}(p)$ is 
$$
-\frac{d}{n-1}\cot(a)+\frac{n-d-1}{n-1}\tan(a)),
$$
and the same is true at all
points of $\tau^{-1}_B(a)$ because of the $G$-invariance. In order to
conclude, we just need to observe that  
$$-\frac{d}{n-1}\cot(a)+\frac{n-d-1}{n-1}\tan(a))\leq 
-\frac{1}{n-1}\cot(a)+\frac{n-1}{n-2}\tan(a)) $$
since $a\in (0,\pi/2)$ and $d\in\{1,\dots,n-1\}$.   
\end{proof}

Reversing the time in the proof of Theorem~\ref{t.barriers-AdS}, one
gets:  

\begin{teo}
\label{t.reverse-barriers-AdS}
Let $E^+(\Lambda)$ be the future tight region of an AdS regular
domain, and $\widehat{\tau}: E^+(\Lambda) \rightarrow (0, \pi/2)$  be
the associated reverse cosmological time. For every $a \in 
(0, \pi/2)$, the generalized mean curvature of the level set
$\widehat{S}_a = \widehat{\tau}^{-1}(a)$ satisfies
$$\frac{1}{n-1}\cot(a) - \frac{n-2}{n-1}\tan(a)\leq
H_{\widehat{S}_a}\leq \cot(a).$$
\end{teo}

\section{CMC time functions in anti-de Sitter spacetimes}
\label{proof.AdS}

The proof follows the same lines as those of Theorem~\ref{t.main-flat-case}
in~\cite[section 6]{ABBZ}), but slightly complicated by the fact that
we need to consider also the reverse cosmological time
(cf. remark~\ref{rk.reverse}). 

\begin{proof}[Proof of Theorem \ref{t.main-AdS-case}]
Let $(M,g)$ be a $n$-dimensional  MGHC spacetimes with constant
curvature $-1$. According to  Theorem~\ref{t.iso-quotient}, $(M,g)$ is
the quotient of a regular domain $E(\Lambda)\subset \mbox{AdS}_n$ by a
torsion-free discrete group $\Gamma\subset\mbox{Isom}(\mbox{AdS}_n)$.
The cosmological time $\tau:E(\Lambda)\rightarrow (0,+\infty)$ and the
reverse cosmological time $\widehat\tau:E(\Lambda)\rightarrow
(0,+\infty)$ are well-defined and regular
(Proposition~\ref{pro.adsregular}).  

For every $a\in [0,+\infty]$, let $S_a=\tau^{-1}(a)$ and $\Sigma_a$ bethe projection  of $S_a$ in $M\equiv\Gamma\setminus E(\Lambda)$. 
Every level set $S_{a}$ is quite obviously a Cauchy hypersurface in
$E(\Lambda)$. It is $\Gamma$-invariant since the cosmological time is
so. It follows that $\Sigma_a$ is a topological Cauchy. hypersurface 
in $M$ since inextendible causal curves in $M$ are projections 
of inextendible causal curves in $E(\Lambda)$. Moreover,
Theorem~\ref{t.barriers-AdS} implies that the generalized
mean curvature $H_{\Sigma_a}$ of $\Sigma_a$ satisfies 
$$-\cot(a)\leq H_{\Sigma_a}\leq
-\frac{\cot(a)}{n-1}+\frac{n-2}{n-1}\tan(a).$$
Consider a decreasing sequence of positive real numbers
$(a_m)_{m\in\NN}$ such that $a_m\to 0$ when $m\to+\infty$. 
Observe that
$$-\frac{\cot(a_m)}{n-1} +
\frac{n-2}{n-1}\tan(a_m)\quad\mathop{\longrightarrow}_{m\to\infty}\quad-\infty.$$
This shows that $(\Sigma_{a_m})_{m\in\NN}$ is a sequence of past
asymptotic $(-\infty)$-barriers.  

For every $a\in [0,+\infty]$, let $\widehat{S}_a=\widehat\tau^{-1}(a)$
and $\widehat{\Sigma}_a$ be the projection  of $\widehat{S}_a$ in
$M$. Of course, $\widehat{\Sigma}_a$ is a topological Cauchy
hypersurface in $M$ for every $a$. By
Theorem~\ref{t.reverse-barriers-AdS}, the generalized
mean curvature $H_{\widehat\Sigma_a}$ of $\widehat\Sigma_a$ satisfies 
$$\frac{1}{n-1}\cot(a) - \frac{n-2}{n-1}\tan(a)\leq H_{\Sigma_a}\leq
\cot(a).$$
Consider a decreasing sequence of positive real numbers
 $(b_m)_{m\in\NN}$ such that $b_m\to 0$ when 
$m\to +\infty$.  Observe that 
$$\frac{1}{n-1}\cot(b_m) - \frac{n-2}{n-1}\tan(b_m)\quad\mathop{\longrightarrow}_{m\to\infty}\quad-\infty.$$
This shows that $(\widehat{\Sigma}_{b_m})_{m\in\NN}$ is a sequence of
past asymptotic $(+\infty)$-barriers. 

So we are in a position to apply Theorem~\ref{t.foliation}, which
shows that $M$ admits a globally defined CMC-time $\tau_{cmc}:M\to
(-\infty,+\infty)$.   
\end{proof}

\section{Description of de Sitter MGHC spacetimes} \label{sec:dS-case} 

We now start our investigation of MGHC de Sitter spacetimes
(i.e. MGHC spacetimes with constant curvature $+1$). Each section in
the sequel is a ``de Sitter substitute'' of a section above dealing with 
anti-de Sitter spacetimes. Our first task will be to introduce a de Sitter
analog of the notion of \emph{AdS regular domain}, called \emph{dS
  standard spacetime}. Every MGHC de Sitter spacetime is
the quotient of a dS standard spacetime by a torsion free subgroup of
$\mbox{Isom}_0(\dS_n)=\mbox{O}_0(1,n)$. Then, we will try to get a good
understanding of the geometry of dS standard spacetimes, 
in order to obtain some estimates of the (generalized) mean curvature
of the level sets of the cosmological time. 

In comparison to the anti-de Sitter case, a major technical
difficulty appears: given a MGHC de Sitter spacetime $(M,g)$, the
developing map $D:\widetilde M\to\dS_n$ is not one-to-one in
general. A consequence is that dS standard spacetimes cannot be
defined as domains in the de Sitter space $\dS_n$. A dS standard
spacetime is a simply connected manifold which is \emph{locally}
isometric  to $\dS_n$; in some particular cases, this manifold is
globally isometric to an open domain in $\dS_n$, but this is not the
general case. 

\subsection{dS standard spacetimes}
\label{subsec.ds}

The purpose of this section is to define a class of locally de Sitter
manifolds, called \emph{dS standard spacetimes}. Recall that a
\emph{M\"{o}bius manifold} is a manifold equipped with a
$(G,X)$-structure, where $X=\SS^{n-1}$ is the $(n-1)$-dimensional sphere
and $G\equiv\mbox{O}_0(1,n)$ is the M\"obius group (i.e. the group of
transformations preserving the usual conformal structure of
$\SS^{n-1}$). To every $(n-1)$-dimensional simply connected M\"obius
manifold $S$, we will associate a $n$-dimensional future complete dS
standard spacetimes $\cB_0^+(S)$ (diffeomorphic to $S\times\RR$). A
similar construction leads to a $n$-dimensional past complete 
dS standard spacetime $\cB_0^-(S)$.

The definition of dS standard spacetimes we will use here first
appeared in a paper by Kulkarni and Pinkall (see \S~3.4 of
\cite{kulkarni}). Unfortunately, Kulkarni-Pinkall did not insist on
the de Sitter nature of the space they consider, and we need to
formulate here the lorentzian interpretation of some of their
results. There is another construction by Scannell (generalizing some
ideas of Mess; see~\cite{scannell} and~\cite{mess}) where the
de Sitter nature of the resulting spaces is obvious. But Scannell
only considered the case of where $S$ is closed, and it is not obvious
from his description that the obtained spacetimes are past or future
complete. So, we will reproduce here Kulkarni-Pinkall's and Scannell's
constructions, for the readers' convenience, and in order to ensure
that both these constructions lead to the same spacetimes.

\subsection{Linear and Klein models of the de Sitter space}
\label{ss.dS}
For $n\geq 2$, let $(x_1,\dots,x_{n+1})$ be the standard coordinate
system on $\RR^{n+1}$, and let $Q_{1,n}$ be the quadratic form $-x_1^2
+ x_1^2 + \ldots + x_{n+1}^2$. The \emph{linear model of the
  $n$-dimensional de Sitter space} is the one-sheeted hyperboloid
$(Q_{1,n}=+1)$ endowed with the lorentzian metric induced by
$Q_{1,n}$; we denote it by $\dS_n$. 

It is easy to check that $\dS_n$ is homeomorphic to
$\RR\times\SS^{n-1}$. Actually, one can prove that $\dS_n$ is
conformally equivalent to
$((-\pi/2,\pi/2)\times\SS^{n-1},-dt^2+ds^2)$, where $dt^2$ is the 
usual metric on $\RR$ and $ds^2$ is the usual metric (of curvature
$1$) on the sphere $\SS^{n-1}$. It follows in particular that $\dS_n$
is globally hyperbolic. The coordinate $x_0$ defines   
on $\mbox{dS}_{n}$ a time function (provided that we make the
appropriate choice of time-orientation).


Observe that each of the two sheets of the hyperboloid $(Q_{1,n}=-1)$
endowed with the riemannian metric induced by $Q_{1,n}$ is a copy of
the $n$-dimensional hyperbolic space. We denote by $H_n^-$
(resp. $H_n^+$) the sheet of the hyperboloid $(Q_{1,n}=-1)$ contained
in the half space $(x_0<0)$ (resp. $x_0>0$).


The projection on $\SS(\RR^{n+1})$ of $\dS_n$ (endowed with the
push-forward of the lorentzian metric of $\dS_n$) is the \emph{Klein
  model of the de Sitter  space}; we denote it by $\DS_n$. The
projections on $\SS(\RR^{n+1})$ of $H^n_-$ and $H^n_+$ will be denoted
by $\HH^n_-$ and $\HH^n_+$. The boundary of $\DS_n$ in
$\SS(\RR^{n-1})$  is the projection of the cone
$(Q_{1,n}=0)\setminus\{0\}$; this is the union of two spheres  
$\SS^{n-1}_+,\SS^{n-1}_-$. We choose the notations such that
$\SS^{n-1}_+$ (resp. $\SS^{n-1}_-$) is included in the projection 
of the half space $x_0>0$ (resp. $x_0<0$). Notice that $\SS^{n-1}_+$
(resp. $\SS^{n-1}_-$) is also the boundary of $\HH^n_-$
(resp. $\HH^n_-$) in $\SS(\RR^{n+1})$.   

Using the conformal structure of $\dS_n$, one sees that every future
oriented inextendible causal curve in $\DS_n$ ``goes from $\SS^{n-1}_-$
to $\SS^{n-1}_+$''. In other words, $\SS^{n-1}_+$ can be seen as the
\textit{future boundary,\/}  of $\DS_n$, and ${\mathbb S}^{n-1}_-$ as
the \textit{past boundary.\/} 

An important observation is that the group $\mbox{O}_0(1,n)$ can be
seen alternatively as the isometry group of the lorentzian space
$\DS_n$, as the isometry group of the 
hyperbolic spaces $\HH^n_-$ and $\HH^n_+$, or as the M\"obius group of
the spheres $\SS^{n-1}_-$ and $\SS^{n-1}_+$ (i.e. the group of the
transformations  preserving the usual conformal structure on the
spheres $\SS^{n-1}_-$ and $\SS^{n-1}_+$). In other words, each
isometry of $\DS_n$ extends as a conformal tranfomation of the spheres
$\SS^{n-1}_-$ and $\SS^{n-1}_+$, and conversely, each conformal
tranformation of the sphere $\SS^{n-1}_\pm$ extends as an isometry of
$\DS_n$. 

\bigskip

The geodesics of $\DS_n$ are the connected components of the
intersections of $\DS_n$ with the projective lines of
$\SS(\RR^{n+1})$. More precisely, let $\gamma$ be a projective line in
$\SS(\RR^{n+1})$, then 
\begin{itemize} 
\item if $\gamma$ does not intersect the spheres $\SS^{n-1}_-$ and
$\SS^{n-1}_+$, then
  $\gamma$ is a spacelike geodesic of $\DS_n$,
\item if $\gamma$ is tangent to the spheres  $\SS^{n-1}_-$ and
$\SS^{n-1}_+$, then each of the two connected components of 
  $\gamma\cap\DS_n$ is a lightlike geodesic in $\DS_n$,
\item if $\gamma$ intersects transversally the spheres 
  $\SS^{n-1}_-$ and $\SS^{n-1}_+$, then each of the two connected components of
  $\gamma\cap\DS_n$ is a timelike geodesic.
\end{itemize}
The causal future $J^+(x)$ of a point $x\in\DS_n$ is
the union of all the projective segments contained in $\DS_n$,
joining at $x$ to $\SS^{n-1}_+$. For the 
timelike future 
$I^+(x)$, one only considers the segments that hit the $\SS^{n-1}_+$
transversally. 
The totally geodesic hypersurfaces in $\DS_n$ are the
connected components of the intersections of $\DS_n$ with the
projective hyperplanes of $\SS(\RR^{n+1})$. 



A key ingredient in the sequel will be the fact that de Sitter space
can be thought of as the space of (non-trivial open) round balls in
${\mathbb S}^{n-1}_+$. For every point $x\in\DS_n$, we denote by
$\partial^+ I^+(x)$ the set of the future endpoints in $\SS^{n-1}_+$
of all the future oriented timelike geodesic rays starting at $x$. 
Then, for every $x\in\DS_n$, the set $\partial^+I^+(x)$ is an open
round ball in $\SS^{n-1}_+$. One can easily check that the map
associating to $x$ the round ball $\partial^+ I^+(x)$ establishes a
one-to-one correspondance between  the points in $\DS_n$ and the
(non-trivial open) round balls in $\SS^{n-1}_+$.
Observe that a point $x\in\DS_n$ is in the (causal) past of another
point $y\in\DS_n$, if and only if the round ball associated to $x$
contains the round ball associated to $y$.
Of course, there is a similar identification between the points of
$\DS_n$ and the round balls in ${\mathbb S}^{n-1}_-$.

\subsection{dS standard spacetimes associated to open domains in
  $\SS^{n-1}_+$}  
\label{ss.particular-case}
Recall that our goal is to associate a future complete dS standard
spacetime $\cB_0^+(S)$ to every simply connected M\"obius manifold 
$S$. In this paragraph, we consider the particular case where $S$ is
an open domain in the sphere $\SS^{n-1}_+$. We denote by $\Lambda$ 
the boundary of $S$ in $\SS^{n-1}_+$. 

For $p\in \Lambda$, let $H(p)$ be the unique projective hyperplane in
$\SS(\RR^{1,n})$ tangent to $\SS_+^{n-1}$ at $p$. Note that
$H(p)\cap\SS^{n-1}_+=\{p\}$, $H(p)\cap\SS^{n-1}_-=\{-p\}$, and $H(p)\setminus \{p, -p \}$ is
contained in $\DS_n$ (more precisely, $H(p)\setminus \{p, -p \}$ is a
lightlike totally geodesic hypersurface in $\DS_n$). Also note that
$\SS(\RR^{n+1})\setminus H(p)$ has two connected components. We denote
by $\Omega^+(p)$ the connected component of $\SS(\RR^{n+1})\setminus
H(p)$ containing $\HH_+^n$. 

\begin{defi}
We consider the set 
$$
\Omega^+(S):=\bigcap_{p \in \Lambda} \Omega^+(p)
$$
We denote by $\cB_0^+(S)$ the unique connected component of
$\Omega^+(S)\cap\DS_{n}$ whose closure in $\SS(\RR^{1,n})$
contains $S$ (see remark~\ref{rk.convex} below). The domain
$\cB_0^+(S)$ is the (future complete) \emph{dS standard spacetime}
associated to $S$. 
\end{defi}

\begin{remark}
\label{rk.convex}
The set $\Omega^+(S)$ is obviously a convex domain of
$S(\RR^{1,n})$. This convex domain contains the hyperbolic space
$\HH^n_+$. Select a point $O\in\HH^n_+$. The radial 
projection of center $O$ on $\SS^{n-1}_+$ 
defines a fibration of $\Omega^+(S)\cap\DS_n$ over
$\SS^{n-1}_+\setminus\Lambda$ with fibers $\RR$. It follows that
there exists a unique connected component of $\Omega^+(S)\cap\DS_{n}$
whose closure contains $S$.  This shows the validity of the above
definition of $\cB_0^+(S)$. 
\end{remark}

\begin{remark}
\label{rk.geod-convex}
Since geodesic segments in $\DS_n$ are segments of
projective lines, another consequence of the convexity of
$\Omega^+(S)$ is the geodesic convexity of $\cB_0^+(S)$:
any geodesic segment joining two elements of $\cB_0^+(S)$ is
contained in $\cB_0^+(S)$.  
\end{remark}

\begin{remark}
\label{rk.causally-convex}
For every $p\in\Lambda$, it is easy to check that the set
$\Omega^+(p)\cap\DS_n$ is the 
timelike 
future of the hyperplane
$H(p)$ in $\DS_n$. It follows that, for every $x\in\cB_0^+(S)$, the
causal future of $x$ in $\DS_n$ is contained in $\cB_0^+(S)$. Since
$\DS_n$ is future complete, it also follows that $\cB_0^+(S)$ is
future complete.  
\end{remark}

\begin{remark}
\label{rk.dual-convex-set}
It is easy to check that, for every $p\in\Lambda$, one has 
$$
\Omega^+(p)=\{ x\in\SS(\RR^{n+1}) \mbox{ such that }\langle x \mid
p\rangle < 0\}
$$
where $\langle \cdot \mid \cdot \rangle$ is the pseudo-scalar product
associated to the quadratic form $Q_{1,n}$. It follows immediately
that $\Omega^+(S)$ is the dual convex set of the convex hull of
$\Lambda$ in $\SS(\RR^{n+1})$. 
\end{remark}

\begin{remark}
\label{rk.domain-dependence}
One can easily check that $\Omega^+(S)\cap\DS_n$ is the set of
points in $\DS_n$ which are not causally related to any element of
$\Lambda$. Therefore, $\cB^+_0(S)$
can be considered as the domain of dependence of $S$ in $\DS_n$, so
that there is a complete analogy between the above definition of dS
standard spacetimes and the definition of AdS regular domains.  
\end{remark}

\begin{remark}
\label{rk.balls}
Recall that there is a canonical identification between the points of
$\DS_n$ and the round balls in $\SS^{n-1}_+$ (see section \ref{ss.dS}). One can
easily check that a point $x\in\DS_n$ is in $\cB_0^+(S)$ if and only
if the ball of $\SS^{n-1}_+$ corresponding to $x$ is contained in
$S$. 
\end{remark}

Of course, there is a similar construction which allow to associate a
past complete domain $\cB_0^-(S)$ to any connected open domain $S$ in
$\SS^{n-1}_-$.

\subsection{The general case}
\label{ss.general-case}
Now we consider the general, where $S$ is any simply connected
$(n-1)$-dimensional M\"obius manifold. A key ingredient will be the
identification between $\DS_n$ and the set of round balls in
$\SS^{n-1}_+$. 

Let us first state two technical lemmas, valid for any local 
homeomorphism $\varphi: X \rightarrow Y$ between manifolds 
(for proofs, see e.g. \cite[\S 2.1]{barHLP}).

\begin{lem}
\label{lem.union}
Let $U$, $U'$ be two open domains in $X$, such that $\varphi$ is
one-to-one in restriction to $U$, and in restriction to $U'$. Assume
that $U \cap U'$  is not empty, and that $\varphi(U')$ contains
$\varphi(U)$. Then, $U'$ contains $U$.\fin
\end{lem}

\begin{lem}
\label{lem.ferme}
Assume that $\varphi$ is one-to-one in restriction to some open domain
$U$ in $X$. Also assume that the set $V = \varphi(U)$ is locally
connected in $Y$, i.e. every point $y$ in the closure of $V$
admits arbitrarly small neighborhood $W$ such that $V \cap W$ is
connected. Then, the restriction of $\varphi$ to the closure of $U$ in
$X$ is one-to-one.\fin 
\end{lem}

Now we start the construction of the dS standard spacetime
$\cB_0^+(S)$. For this purpose, we choose a
$d:S\to\SS^{n-1}_+$. Recall that such a map does exist since $S$ is a
M\"obius manifold. Also recall that the map $d$ is not one-to-one in
general. 

\begin{defi}
An (open) \emph{ round ball} $U$ in $S$ is an open domain in $S$ such
that the developing map $d$ to $U$ is one-to-one in restriction to
$U$, and such that $d(U)$ is an open round ball in $\SS^{n-1}_+$. A
round ball $U\subset S$ is said to be \emph{proper} if the image under
$d$ of the closure $\overline{U}$ of $U$ in $S$ is the closure of
$d(U)$ in $\SS^{n-1}_+$.  
\end{defi}

Note that according to Lemma \ref{lem.ferme}, if $U$ is a proper 
round ball in $S$, then $d$ is one-to-one in restriction to
$\overline{U}$ and $d(\overline{U})$ is a closed round ball of
${\mathbb S}^{n-1}$.  

\begin{defi}
We will denote by $\cB(S)$ the set of all round balls in $S$, and by
$\cB_0(S)$ the set of proper round balls. 

The sets $\cB(S)$ and $\cB_0(S)$ are naturally ordered by the
inclusion. For every element $U$ of $\cB_0(S)$, we denote by $W(U)$
the subset of $\cB_0(S)$ made of the proper round balls $U'$ such that
$\overline{U'}\subset U$. Given two elements $U,V$ of $\cB_0(S)$ such
that $\overline{U}\subset V$, we denote by $W(U,V)$ the set of all
proper round balls $U'$ in $S$ such that $\overline{U}\subset U'$ and
$\overline{U'}\subset V$. The sets $W(U,V)$ generate a topology on
$\cB_0(S)$ that we call the \textit{Alexandrov topology}.
\end{defi}

We already observed that the de Sitter space $\DS_n$, as a set, is
canonically identified with the space ${\mathcal B}_0({\mathbb
  S}^{n-1}_+)={\mathcal B}({\mathbb S}^{n-1}_+)$ of all open round balls
in the sphere $\SS^{n-1}_+$ (see \S\ref{ss.dS}).

\begin{lemma}
\label{l.homeo}
The canonical identification between $\DS_n$ and ${\mathcal B}_0({\mathbb
  S}^{n-1}_+)$ is an homeomorphim, once ${\mathcal B}_{0}({\mathbb
  S}^{n-1}_+)$ is endowed with the Alexandrov topology.
\end{lemma}

\begin{proof}
Let $U,V$ be two points elements in $\cB_0(\SS^{n-1}_+)$ such that
$U\subset V$. Let $x,y$ be the points of $\DS_n$ corresponding
respectively to $U$ and $V$. Recall that this means that $U$
(resp. $V$) is the set of future extremities of timelike geodesics
starting at $x$ (resp. $y$). Hence $U\subset V$ implies $J^+(x)\subset
I^+(y)$, or equivalently $p\in I^+(q)$. Now observe that the set
$W(U,V)\subset \cB_0(\SS^{n-1}_+)$ corresponds in $\DS_n$ to the set
of all points $z$ such that $J^+(x)\subset I^+(z)$ and $J^+(z)\subset
I^+(y)$, or equivalently, $z\in I^+(y)\cap I^-(x)$. But since $\DS_n$
is strongly causal, the topology on $\DS_n$ generated by sets of the
type $I^+(y)\cap I^-(x)$ is the same as the manifold topology. The
lemma follows.
\end{proof}


\begin{prop}
\label{pro.Bmanifold}
The set $\cB_0(S)$, equipped with the Alexandrov topology, is a manifold.
\end{prop}

\begin{proof}[Sketch of proof]
Compare our proof with \cite[Proposition page $98$, item
(iii)]{kulkarni}. The developing map $d:S\to\SS^{n-1}$ induces a map 
$d:\cB_0(S)\to\cB_0(\SS^{n-1})$. The composition of this map with the
identification between $\cB_0(\SS^{n-1}_+)$ with $\DS_n$ defines a
natural map $\cD^+:\cB_0(S)\rightarrow\DS_n$.  
For any element $U$ of $\cB_0(S)$, the restriction of $F$ to $W(U)$ 
is a homeomorphism onto its image, which is the future $I^+(x)$ of the
point $x$ such that  $\partial{I}^+(x) = d(U)$. It follows that the
$W(U)$ are charts on $\cB_0(S)$ homeomorphic to $\RR^{n}$. 

Let us prove the Hausdorff separation property: let $U_1$, $U_2$ be
elements of $\cB_0(S)$ such that every neighborhood of $U_1$
intersects every neighborhood of $U_2$. Let $U'_1$, $U'_2$ be other
elements of $\cB_0(S)$ such that  $\overline{U}_1 \subset U'_1$
and $\overline{U}_2 \subset U'_2$. Then, the neighborhoods $W(U'_1)$
and $W(U'_2)$ have non-trivial intersection since the first contains
${U}_1$ and the second contains ${U}_2$. Let $V$ be a common
element. The round ball $V$ is contained in $U'_1 \cap U'_2$. This
last intersection is not empty. According to  Lemma \ref{lem.union},
the image by $D$ of $U'_1 \cap U'_2$ is $D(U'_1) \cap D(U'_2)$. It follows 
that the restriction of $D$ to the union $U'_1 \cup U'_2$ is injective. 
Therefore, the restriction of $\cD^+$ to $W(U'_1) \cup W(U'_2)$ is a
homeomorphism, and $\cD^+(W(U'_1)\cup W(U'_2))=I^+(\cD^+(U'_1)) \cup
I^+(\cD^+(U'_2))$. Since the Hausdorff property holds in
$I^+(\cD^+(U'_1)) \cup I^+(\cD^+(U'_2))$, we conclude that $U_1 = U_2$. 

The fact that $\cB_0(S)$ is second countable is not really relevant to
our purpose, and its proof is left to the reader.
\end{proof}

The map $\cD^+:\cB_0(S) \rightarrow\DS_n$ (obtained as the composition
of the developing map $d:\cB_0(S)\rightarrow\cB_0(\SS^{n-1})$ and the 
identification of $\cB_0(\SS^{n-1})$ with $\DS_n$) is a local
homeomorphism (see Lemma~\ref{l.homeo}). Hence, we can consider the
pull-back by $\cD^+$ of the de Sitter metric on $\cB_0(S)$. This is a
locally de Sitter lorentzian metric on $\cB_0(S)$. 


\begin{defi}
\label{d.dS-regular-domain}
We will denote by $\cB_0^+(S)$ the manifold
$\cB_0(S)$ equipped with the pull-back by $\cD^+$ of the de Sitter metric.
\end{defi}

\begin{remark}
It is clear from our definitions that the lorentzian manifold
$\cB_0^+(S)$ is future complete. It is also that $\cB_0^+(S)$ is 
\emph{asymptotically simple,\/} i.e. that every inextendible
future oriented null geodesic ray is complete. It follows
that $\cB_0^\pm(S)$ is globally hyperbolic  (see Proposition $2.1$ in
\cite{dscft}). 
\end{remark}

\begin{prop}
\label{prop.def-coincide}
In the case where the map $d$ is one-to-one, the lorentzian manifold
$\cB_0^+(S)$ defined in this paragraph is isometric to the domain
$\cB_0^+(d(S))$ defined in \S~\ref{ss.particular-case}.    
\end{prop}

\begin{proof}
This follows immediately from the constructions and from
remark~\ref{rk.balls}. The isometry is given by the map $\cD^+$. 
\end{proof}

\begin{remark}
If $d':S\to\SS^{n-1}_+$ is another developing map, then $d'=\phi\circ
d$ where $\phi$ is an element of the M\"obius group $O(1,n-1)$ (in
particular, $\phi$ maps round balls on round balls). It follows that, 
up to isometry, the dS standard spacetime $\cB_0^+(S)$ does depend on
the choice of $d$.
\end{remark}

A similar construction (where the sphere $\SS^{n-1}_+$ is replaced by
the sphere $\SS^{n-1}_-$) yields a past complete lorentzian manifold
$\cB_0^-(S)$.

\begin{defi}
A future (resp. past) complete \emph{dS standard spacetime}
is a lorentzian manifold of the type $\cB_0^+(S)$ (resp. $\cB_0^-(S)$)
where $S$ is a simply connected M\"obius manifold.   
\end{defi}

If $S$ is conformally equivalent to the sphere ${\mathbb S}^n$, then
$S$ and $\cB_0^\pm(S)$ is said to be \textit{elliptic.\/} If $S$ is
conformally equivalent to the sphere ${\mathbb S}^n$ minus a single
point,  then $S$ and $\cB_0^\pm(S)$ are said to be
\textit{parabolic}. If $S$ is neither elliptic nor parabolic, then $S$
and $\cB_0^\pm(S)$ are said to be \textit{hyperbolic}.

\begin{remark}
\label{rk.elementary}
According to these definitions, there is only one elliptic standard 
dS spacetime: the de Sitter space itself. Up to isometry, there is only 
one future complete (resp. past complete) parabolic standard
spacetime, which can be described as the future (resp. past) in
$\DS_n$ of a point in the conformal boundary $\SS^{n-1}_-$
(resp. $\SS^{n-1}_+$).
\end{remark}

\subsection{Canonical neighbourhood and canonical domain of a point}
Let $S$ a simply connected M\"obious manifold of dimension $n-1$. Let
$d:S\to\SS^{n-1}_+$ be a developing map. In general, the dS
standard spacetime $\cB_0^+(S)$ does not admit any global isometric
embedding in $\DS_n$. Nevertheless, for many purpose, we will not need
to study the geometry of the whole spacetime $\cB_0^+(S)$,
but only the geometry of some regions of $\cB_0^+(S)$ (typically
the past of a point in $\cB_0^+(S)$). The purpose of this paragraph is
to define some ``big'' regions of $\cB_0^+(S)$ which admit some
isometric embeddings in $\DS_n$.  

\begin{defi}
For $x\in S$, we denote by $\cU(x)$ the union of all the open round
balls containing $x$. The set $\cU(x)$ is called the \emph{canonical
  neighborhood} of $x$ in $S$. 
\end{defi}

Using Lemma~\ref{lem.union} it is easy to prove the following
proposition (see also \cite[Proposition~4.1]{kulkarni}).

\begin{prop}
\label{pro.canonicalneigh}
The restriction of $d$ to any canonical neighborhood is one-to-one. 
\fin
\end{prop}




Putting together Propositions~\ref{pro.canonicalneigh}
and~\ref{prop.def-coincide}, we get.

\begin{coro}
\label{c.canonical-domain}
For every $x\in S$, the dS standard spacetime $\cB_0^+(\cU(x))$ is
isometric to the dS standard spacetime $\cB_0^+(d(\cU(x)))$
(associated to the open domain $d(\cU(x))$ of $\SS^{n-1}_+$. In
particular, $\cB_0^+(\cU(x))$ is globally isometric to an opain domain
in $\DS_n$. 
\end{coro}

Moreover, the past of a point can always be seen in a domain of the
form $\cB_0^+(\cU(x))$.
 
\begin{prop}
\label{pro.voishorizon}
Let $U$ be an element of $\cB^+_0(S)$ (i.e. a proper round ball in
$S$). Let $x\in U\subset S$. Then the canonical domain
$\cB_0^+(\cU(x))$ contains the past of $U$ in $\cB_0^+(S)$. 
\end{prop}

\begin{proof}
Recall that a round ball $V$ is in the past of $U$ in $\cB_0^+(S)$
if and only if $V$ contains $U$. So, if $V$ is in the past of
$U$, then $x\in V$; hence, $V\in\cB_0^+(cU(x))$. 
\end{proof}

\subsection{Another definition of dS standard spacetimes}
The construction of dS standard spacetimes detailed in the previous
paragraph is quite different from those given by Scannell
in~\cite{scannell}. We will now explain Scannell's constuction. 

\medskip

Let $S$ be a \emph{hyperbolic} simply connected M\"obius manifold of
dimension $n-1$, and $d:S\to\SS^{n-1}_+$ be a developing map. Let
$\cB_{max}(S)$ be the set of maximal open round balls in $S$, i.e. the
maximal elements of $\cB(S)$.  
For every element $U$ of $\cB_{max}(S)$, let $\overline{U}$ be the 
the closure of $U$ in $S$, let $\overline{d(U)}$ be the closure of
$d(U)$ in $\SS^{n-1}$, and let $\Lambda_S(U)$ be the complement of
$d(\overline{U})$ in  $\overline{d(U)}$. Observe that $\Lambda_S(U)$
is closed in ${\mathbb S}^{n-1}$. The closed set $\overline{d(U)}$ is
conformally equivalent to the compactified hyperbolic space
$\HH^{n-1}\cup\partial\HH^{n-1}$. We may therefore transfer the usual
notion of hyperbolic convex hull to $\overline{d(U)}$, and define the 
convex hull $\hat{C}(U)$ of $\Lambda_S(U)$ in $\overline{d(U)}$. Let
$C(U)=d^{-1}(\widehat{C}(U))\cap U$ (note that $C(U) =
\emptyset$ if and only if $\Lambda_S(U)$ has less  than two points).
A key point in the construction is the following fact 
(\cite[Theorem 4.4]{kulkarni} or \cite[Proposition~4.1]{scannell}). 

\begin{fact}
For every $x$ in $S$ there exists a unique element $U(x)$
of $\cB_{max}(S)$ such that $x$ belongs to $C(U(x))$.  
\end{fact} 

\begin{rema}
This fact allows to define a stratification of the M\"obius manifold $S$: for every $x\in S$, the strata of $x$ is the set $C(U(x))$. This stratification --- which was defined by Thurston in some particular case (unpublished), and later by Apanasov and Kulkarni-Pinkall in the general case (\cite{kulkarni}) --- is called the \emph{canonical stratification} of $S$. 
\end{rema}

Following Scannell (see~\cite[page 8]{scannell}), we will now define a 
local homeomorphism ${\bf D}^+: S \times (0, +\infty ) \rightarrow
\DS_n$. We use the identification the points in $\DS_n$ and the set of
round balls in $\SS^{n-1}_+$: for every $x$ in $S$, we see the round ball
$U(x)$ as a point in $\DS_{n}$. Let $c_x:
[0,+\infty)\rightarrow\DS_n$ 
be the unique unit speed future oriented timelike geodesic such that
$c(0)=U(x)$ and $c(t) \to x$ when $t \to \infty$. We define ${\bf
  D}^+(x, t)$ as the point $c_x(t)$ in $\DS_n$. Scannell proved that this
map is a local homeomorphism. Then we can define the future complete
\emph{dS standard spacetime} $\mathbf{B}^+(S)$ associated with $S$ as the
manifold $S\times (0, +\infty)$ equipped with the pull-back by ${\bf
D}^+$ of the de Sitter metric. 

We will see later (remark \ref{rk.memesca}) that this definition of
dS standard spacetimes coincides with the definition given in
\S\ref{ss.general-case} (more precisely, the locally de Sitter
manifolds $\mathbf{B}^+(S)$ and $\cB_0^+(S)$ are isometric). At this
point, it should be clear to the reader that there exists an isometric
embedding $f: {\bf B}^+(S) \hookrightarrow {\mathcal B}^+_0(S)$ such that 
${\bf D}^+ = {\mathcal D}^+ \circ f$. 

\subsection{MGHC de Sitter spacetimes and dS standard spacetimes}

The reason of being of dS standard spacetimes is the following Theorem: 

\begin{teo}[Scannell]
\label{teo.dscompact}
Every  MGHC  dS-spacetime is the quotient of a dS standard spacetime
by a torsion-free discrete subgroup of isometries.
\end{teo}

\begin{proof} 
See~\cite{scannell} (and remark~\ref{rk.memesca} which shows that Scannell's definition of dS standard spacetimes is equivalent to Kulkarni-Pinkall's definition).
\end{proof}

\section{Cosmological time and horizons of dS standard spacetimes}
\label{s.singularity}

All along this section, we consider a simply connected M\"obius
manifold $S$ of dimension $n-1$, and the associated (future complete)
dS standard spacetime $\cB_0^+(S)$. We assume that $S$ is hyperbolic. 

Recall that $\cB_0^+(S)$ is defined as follows. One chooses  a
  developing map $d:S\to\SS^{n-1}_+\simeq\SS^{n-1}$. One considers the
  space $\cB_0(\SS^{n-1}_+)$. This map induces a local homeomorphism
  $d:\cB_0(S)\to\cB_0(\SS^{n-1}_+)$. The composition of this local
  homeomorphism  with the identification between $\DS_n$ with
  $\cB_0(\SS^{n-1}_+)$ defines a local homeomorphism
  $\cD^+:\cB_0(S)\to\DS_n$. The dS standard spacetime 
  $\cB_0^+(S)$ is, by definition, the manifold $\cB_0(S)$ equipped 
  with the pull back by $\cD^+$ of the lorentzian metric of $\DS_n$.
So, by construction, $\cD^+$ defines a locally isometric developing
  map of $\cB_0^+(S)$ in $\DS_n$. 

The purpose of this section is to get some informations on the
cosmological time of $\cB_0^+(S)$. Just as in the AdS setting, this
will lead us to study the support hyperplanes of the past horizon
$\cH^-(S)$ of $\cB_0^+(S)$. Of course, a similar study could be
carried out for the dS standard spacetime $\cB_0^-(S)$.

\subsection{Cosmological time}

\begin{prop}
\label{pro.dsregular}
The dS standard spacetime $\cB_0^+(S)$ has a regular cosmological
time.  
\end{prop}

\begin{proof}
Recall that we have assumed that $S$ is hyperbolic; this will
play a crucial role here. We denote by $\tau$ the cosmological time
of $\cB_0^+(S)$.    

Let $x\in\cB_0^+(S)$. We want to prove that $\tau(x)$ is finite. We
argue by contradiction. If $\tau(x)=+\infty$, then, for every
$n\in\NN$, we can find a past directed causal curve
$c_n:[0,1]\to\cB_0^+(S)$ such that $c_n(0)=x$ and such that the length
of $c_n$ is at least $n$. For every $n$, let $x_n:=c_n(1)$. Let
$z:=\cD^+(x)$. For every $n\in\NN$, let $\gamma_n:=\cD^+\circ c_n$ and
$z_n:=\cD^+(x_n)=\gamma_n(1)$. Then $(\gamma_n)_{n\in\NN}$ is a 
sequence of past directed compact causal curves in $\DS_n$, all having
the same final extremity $z$, and such that the length of $\gamma_n$
tends to $\infty$ when $n\to\infty$. It follows that, up to extracting
a subsequence, the sequence $(z_n)_{n\in\NN}$ converges to a point
$\bar{z}\in\SS^{n-1}_-$. Now, recall that $(x_n)_{n\in\NN}$ is a
sequence of points in $\cB_0^+(S)$, that is, a sequence of proper 
round balls in $S$. Let $\bar x$ be the liminf of these balls,
\emph{i.e.} $\bar x=\bigcup_{p\in\NN}\bigcap_{n\geq p} x_n.$ 
Note that $d$ is one-to-one in restriction to $\bar x$ (since it is
one-to-one in restriction to each $x_i$). For every $n$, the point
$z_n$ can be seen as a ball in $\SS^{n-1}_+$ 
(using the identification of $\DS_n$ with the space of round balls in
$\SS^{n-1}_+$). If we see $x_n$ as a ball in $S$ and 
$z_n$ as a ball in $\SS^{n-1}_+$, then we have $z_n=d(x_n)$. Hence,
$d(\bar x)$ is the liminf of the sequence of balls
$(z_n)_{n\in\NN}$. Since $z_n\to\bar z\in\SS^{n-1}_-$, it follows that
$d(\bar x)$ is the complement of a single point in
$\SS^{n-1}_+$. According to Lemma \ref{lem.ferme},  
this implies that the boundary of the ball $\bar x$ in $S$ is either
empty, or a single point. In the former case, we have $y=S$, hence $S$ 
is parabolic, and this contradicts our hypothesis on $S$. In the
latter case, the restriction of $d$ to the closure $\overline{y}$ is a
homeomorphism onto $\SS^{n-1}_+$; it follows that $S$ is
elliptic, and this also contradicts our hypothesis. So we have that
$\tau(x)$ is finite. 

Now, we consider an inextendible past oriented causal curve
$c:[0,T)\to\cB_0^+(S)$. We have to prove that $\tau(c(t))\to 0$ when  
$t\to T$. Let $x:=c(0)$. On the one hand, for every $t\in [0,T)$, the
quantity $\tau(c(t))$ does not depend on the whole spacetime
$\cB_0^+(S)$, but only on the past $J^-(x)$ of $x$ in $\cB_0^+(S)$. On
the other hand, the set $J^-(x)$ is contained in the domain
$\cB_0^+(\cU(x))$  (Proposition~\ref{pro.voishorizon}). As a
consequence, in our problem, we can replace the cosmological time
$\tau$ of the dS standard spacetime $\cB_0^+(S)$ by the cosmological
time $\check\tau$ of standard spacetime
$\cB_0^+(\cU(x))\subset\cB_0^+(S)$. But the standard spacetime
$\cB_0^+(\cU(x))$ is isometric to a causally convex domain of $\DS_n$
(corollary~\ref{c.canonical-domain} and
remark~\ref{rk.causally-convex}). It follows easily that 
$\check\tau(c(t))\to 0$ when $t\to T$. Therefore $\tau(c(t))\to 0$
when $t\to T$.   
\end{proof}

\begin{remark} 
\label{rk.cosmoele}
\begin{enumerate} 
\item 
Since MGHC de Sitter spacetimes are quotients of standard spacetimes 
by Theorem \ref{teo.dscompact}, 
and since cosmological time functions are preserved by 
isometries, it is an immediate corollary of Proposition \ref{pro.dsregular} 
that MGHC hyperbolic standard spacetimes have regular cosmological time.

In \cite[Theorem 3.1]{dscft} it is shown that for a class 
of MGHC spacetimes (spacetimes of de Sitter type), satisfying the
strong energy condition with positive
cosmological constant, assuming that the future conformal boundary has an
infinite fundamental group implies that the spacetime is past incomplete. 

This result and our Proposition \ref{pro.dsregular} 
have quite similar flavor. 
The result in \cite{dscft} is more general since 
MGHC spacetimes of de Sitter type do not have in general constant curvature. 
On the other hand, the conclusion of Proposition 8.1 is stronger, since a
spacetime may be past incomplete without having a regular cosmological time. 

\item Elliptic and parabolic dS standard spacetimes do not have regular
cosmological time. The cosmological time in these spacetimes is
everywhere infinite.  
\end{enumerate}
\end{remark}

Of course, there are analogs of Theorem \ref{teo.cosmogood} and
Proposition \ref{pro.dsregular}, concerning the reverse cosmological
time in past complete dS standard spacetimes.

\subsection{Past horizon}

As in the AdS case, one can define a notion of \emph{past horizon} for
future complete dS standard spacetimes. Recall that $\cB_0(S)$ is the
set of \emph{proper} open round balls in $S$, whereas $\cB(S)$ is the
set of \emph{all} round balls in $S$ (see \S~\ref{ss.general-case}).

\begin{defi}
The \emph{past horizon} of the future complete regular domain
$\cB_0(S)$ is the set $\cH^-(S):=\cB(S)\setminus\cB_0(S)$. 
\end{defi}

\begin{rema}
\begin{enumerate}
\item  The arguments of Proposition~\ref{pro.Bmanifold} can be easily
  adapted, leading to the conclusion that the set $\cB(S)$ admits a
  topology for which it is a manifold with boundary (the boundary
  being precisely the past horizon
  $\cH^-(S)=\cB(S)\setminus\cB_0(S)$). Moreover, the developing map
  $\cD^+:\cB_0(S)\to\DS_n$ extends to a local homeomorphism from
  $\cB(S)$ into $\DS_n$ that we still denote by $\cD^\pm$.
\item Every round ball in $S$ is the increasing union of one-parameter
  family of proper round balls. It follows that any past-extendible
  causal curve $c$ in $\cB_0^+(S)$ admits a limit point in  the
  horizon $\cH^-(S)$; we call this point the \emph{initial 
  extremity} of the curve $c$. Conversely, any point $p\in\cH^-(S)$ is
  the initial extremity of a past-inextendible timelike curve in
  $\cB_0^-(S)$ (which can actually be chosen to be geodesic). 
\item Recall that, in the particular case where $S$ is an open domain
  in $\SS^{n-1}_+$, the dS standard spacetime $\cB_0^+(S)$ can be seen
  as an open domain domain in the de Sitter space $\DS_n$. Using
  item~(1), it is easy to see that, in this 
  particular case, the past horizon $\cH^-(S)$ is just the topological
  boundary in $\DS_n$ of the open domain $\cB_0^+(S)$. 
\end{enumerate}
\end{rema}



As noticed above, the past horizon $\cH^-(S)$ admits a simple
description in the particular case where the developing map $d$ is
one-to-one. Lemma~\ref{lem.U} shows that, as far as ``semi-local''
properties of $\cH^-(S)$, one can always reduce to this particular
case. We recall that every point $q\in S$ admits a ``nice''
neighbourhood $\cU(q)$ in $S$ which is isometric to an open domain in
$\SS^{n-1}_+$. 

\begin{lem}
\label{lem.U}
Let  $p$ be a point in $\cH^-(S)$. Let $c$ be a future complete
timelike geodesic with initial extremity $p$. Let $q$ be the future
extremity of $c$ in ${\mathbb S}^+_{n-1}$. For every element $x$ in
$c$, let $\cH^-_x(S)$ be the intersection of $\cH^-(S)$ with the
closure of $I^-(x)$ in $\cB^+(S)$. Similarly, let $\cH^-_y(\cU(q))$ be 
the intersection of $\cH^-(\cU(q))$ with the closure of $I^-(x)$ in
$\cB^+(\cU(q))$. Then $\cH^-_x(S)$ is an open neighborhood of $p$ in 
$\cH^-(S)$ and coincides with $\cH^-_x(\cU(q))$.
\end{lem}

\begin{proof}
This is an immediate corollary of  Proposition \ref{pro.voishorizon}.
\end{proof}


Let us assume that $S$ is a domain in the sphere $\SS^{n-1}_+$. Recall
that, under this assumption, the dS standard spacetime 
$\cB_0^+(S)$ is a domain in $\DS_n$, and the past horizon $\cH^-(S)$
is just the boundary of $\cB_0^+(S)$ in $\DS_n$. Also recall that
$\cB_0^+(S)$ is defined as a connected component of the intersection
of the convex set $\Omega^+(S)$ with $\DS_n$ (see \S~\ref{ss.particular-case}). In
particular, $\cH^-(S)$ is a locally convex hypersurface in
$\SS(\RR^{n+1})$. This allows us to speak of the support 
planes of $\cH^-(S)$ (which are projective hyperplanes in
$\SS(\RR^{n+1}$). Note that, just as in AdS case, if $H$ is a support
hyperplane of $\cH^-(S)$, then the totally geodesic hypersurface
$H\cap\DS_n$ is a spacelike or degenerate. The following statement is
the analog of Proposition~\ref{lem.adshori} in the AdS case.

\begin{prop}
\label{lem.dshori}
Assume that $S$ is a domain in $\SS^{n-1}_+$. Let $p$ a point of
${\mathcal H}^-(S)$. Let $C(p)\subset T_p\DS_n$ be 
the set of the future directed unit tangent vectors orthogonal to the
support hyperplanes of $\cH^-(S)$ at $p$. Then:  
\begin{enumerate}
\item the set $C(p)$ is the convex hull of its lightlike elements;
\item If $c$ is a future complete geodesic ray starting at $p$ whose
  tangent vector at $p$ is a lightlike element of $C(p)$, then the
  future endpoint of $c$ is in $\Lambda$ (recall that $\Lambda$ is the
  boundary of $S$ in $\SS^{n-1}_+$). 
\end{enumerate}
\end{prop}

\begin{proof}
The proof is very similar to those of Proposition~\ref{lem.adshori};
the only differences are the following.
\begin{itemize} 
\item We work with the convex set $\Omega^+(S)$ instead of the convex
  set $E(\Lambda)$.
\item The point $q$ now belongs to $\HH^n_+\cup\SS^{n-1}_+$ (instead
  of $\ADS_n\cup\partial\ADS_n$ in the AdS case).
\item  The causal vector $v_q$ is lighlike if and only if
$q\in\SS^{n+1}_+$.
\item The proof of item~(2) is slightly easier in the dS case: since
  the quadratic form $Q_{1,n}$ has signature $(1,n)$, one gets that
  the subspace spanned by the $\hat q_i$'s is 1-dimensional (instead
  of $2$-dimensional in the AdS case); it follows immediately that all
  the $q_i$'s are equal to $q$, and thus, that $q$ is in $\Lambda$. 
\end{itemize} 
\end{proof}

\subsection{Retraction onto the horizon}


We will now study the realizing geodesics in $\cB_0^+(S)$.
Let $x\in\cB_0^+(S)$. Recall that a future directed timelike geodesic ray
$c:(0,1]\to\cB_0^+(S)$ such that $c(1)=x$ is a \emph{realizing
  geodesic for $x$} if $\tau(x)$ is equal to the length  
of $c$. Clearly, realizing geodesic rays for $x$ are contained in
the past of $x$. Therefore, for our problem, we may pick a point
$q\in\SS^{n-1}_+$ which is the future endpoint of a timelike geodesic
passing through $x$, and replace the dS standard spacetime
$\cB_0^+(S)$ by the dS standard spacetime $\cB_0^+(\cU(q))$
(Proposition~\ref{pro.voishorizon}). In other words, as far as
realizing geodesic rays for $x$ are concerned, we may assume without
loss of generality that $S$ is an open domain in the sphere $\SS^{n-1}_+$.   

\begin{prop}
\label{pro.dsunique}
For every $x\in\cB_0^+(S)$, there is a unique realizing geodesic
for $x$ in $\cB_0^+(S)$.   
\end{prop}

\begin{proof}
Recall that we assume (without loss of generality) that $S$ is a
domain in $\SS^{n-1}_+$. Hence, the dS standard 
spacetime $\cB_0^+(S)$ is a connected component of the intersection of
the convex set $\Omega^+(S)$ with $\DS_n$, and $\cH^-(S)$ is the
boundary of $\cB_0^+(S)$ in $\DS_n$.  
Initial extremities of realizing geodesics for $x$ are points $z$ in 
$\cH^-(S)$ such that $d(x,z)=\tau(x)$, where
$d(x,z)$ is the length of a past oriented timelike geodesic in $\DS_n$
starting from $x$ and ending to $z$. For each $\tau$, the set 
$\{ z \in \DS_n  | d(x,z) \geq \tau \}$ 
is the intersection of  $\DS_n$  with a solid ellipsoid $\cE_\tau$ in
$\SS(\RR^{n})$ tangent to the sphere $\SS_-^{n-1}$ 
along a round subsphere.  If $\tau < \tau'$, 
then $\cE_{\tau'} \subset \mbox{int}\cE_\tau$,  leading to 
the  definition:
$$
\tau(x) = \sup \{ \tau  | \cE_\tau \cap \cH_-(S) \neq 
\emptyset \}.
$$
Let $y$, $y'$ be initial extremities of realizing geodesics for $x$:
they both belong  to $\cE_{\tau(x)} \cap
\overline{\Omega}^+(\Lambda)$. On one hand, the segment $[y, y']$ is
contained in  the interior of ${\mathcal E}_{\tau(x)}$ (since
ellipsoids are strictly convex). On the other hand, according to
Remark \ref{rk.convex}, the segment $[y,y']$  is contained in
${\mathcal B}^+(S)$. We obtain a contradiction, unless $y = y'$ (see
the proof of Proposition~\ref{pro.tout}). 
\end{proof}

\begin{prop}
\label{pro.qdrealise}
Let $c: (0, T] \rightarrow \cB_0^-(S)$ be a future oriented timelike
geodesic whose initial extremity $p=\lim_{t\to 0}c(t)$ belongs to the
past  horizon $\cH^-(\Lambda)$. Then the following assertions are
equivalent.  
\begin{enumerate}
\item the geodesic $c$ is tight,
\item there exists $t_0\in (0, T]$ such that $c((0,t_0])$ is a realizing
  geodesic for the point $c(t)$,  
\item $c$ is orthogonal to a support hyperplane of $\Omega^+(S)$ at
  $p$.  
\end{enumerate}
\end{prop}


\begin{proof}
The proof is entirely similar to those of Proposition \ref{pro.tight}, 
based on the strict convexity of the ellipsoids $\cE_\tau$.
\end{proof}



\begin{remark}
\label{rk.memesca}
According to Lemma~\ref{lem.dshori}, and since there is at least one
realizing geodesic for each $x$ in ${\mathcal B}^+_0(S)$,
Proposition~\ref{pro.qdrealise} means precisely  that the map  
$f: {\bf B}(S) \rightarrow\cB^+_0(S)$ defined at the end of 
\S~\ref{subsec.ds} is onto. Hence $f$ is an isometric identification
between $\cB^+_0(S)$ and ${\bf B}(S)$. 
\end{remark}

\section{Curvature estimates of cosmological levels in dS standard spacetimes}

\begin{teo}
\label{t.barriers-dS}
Let ${\mathcal B}^+_0(S)$ be a future complete dS standard spacetime,
and $\tau: {\mathcal B}^+_0(S) \rightarrow (0, +\infty)$  
be the associated cosmological time function. Then, for every $a\in
(0,+\infty)$, the generalized mean curvature of the level set $S_a =
\tau^{-1}(a)$ admits the following estimates
$$
-\coth(a)\leq H_{S_{a}} \leq -\frac{1}{n-1}\coth(a) -
\frac{n-2}{n-1}\tanh(a).
$$ 
\end{teo}

\begin{proof}
We use the same notations $x$, $p$, $c$, $v$ as in the proof of
Theorem \ref{t.barriers-AdS}. The past of the geodesic
$c:\RR\to\cB_0^+(S)$ contains the past in $\cB_0^+(S)$ of a small 
neighbourhood $U$ of $x$. The restriction to $U$ of the function
$\tau$ only depends of the past of $U$ in $\cB_0^+(S)$. Hence the
geometry of the hypersurface $S_a$ in $U$ (in particular the
generalized mean curvature of $S_a$ at $p$) only depends on the past
of $c$ in $\cB_0^+(S)$. Together with Lemma~\ref{lem.U}, this allows
us to restrict ourselves  to the case where $S$ is an open domain in
$\SS^{n-1}_+$. 

The proof is then formally completely similar to those of
Theorem~\ref{t.barriers-AdS}. The hypersurface ${\mathbb S}_x^+$ is
the set of the points of $\cB_0^+(S)$ which are in the future of $p$,
at distance exactly $a$ from $p$. Clearly, $\SS_x^+$ is in the future
of $S_a$, and $x\in\SS_x^+$. A simple computation shows that the mean
curvature of $\SS_x^+$ is constant and equal to $-\coth(a)$.

In order to construct the hypersurface ${\mathbb S}^-_x$, we select a
finite set $v_1,\dots,v_r$ of lightlike elements of $C(p)$ such that
$v\in\mbox{Conv}(v_1,\dots,v_r)$ (such a finite
set does exist by item~(1) of Proposition~\ref{lem.dshori}). For every $i$, we
denote by $q_i$ the future endpoint of the lightlike geodesic ray
whose tangent vector at $p$ is the vector $v_i$. Let
$S'=\SS^{n-1}_+\setminus\{q_1,\dots,q_r\}$. Item~(2)
of Proposition~\ref{lem.dshori} shows that $S'\supset S$. The domain  
$\cB^+_0(S')\subset\DS_n$ is a dS standard spacetime 
with regular cosmological time $\tau'$. We define the hypersurface
$\SS_x^-$ as the $a$-level of the cosmological time $\tau'$. Since
$S'\supset S$, the domain $\cB_0^+(S')$ contains the domain
$\cB_0^+(S)$, and thus, $\SS_x^-$ is in the past of $S_a$. 
So, we are left to compute the mean curvature of $\SS_x^-$ at $p$.
For this purpose, we introduce the minimal projective subspace $F$ in
  $\SS(\RR^{n+1})$ containing $q_1,\dots,q_r$. We observe that
  $\SS_x^-=(\tau')^{-1}(a)$ 
is the saturation under $G$ of the umbilical submanifold $S_a\cap
F^\perp$, where $G$ is the group of isometries fixing $F$ 
pointwise. It follows that the mean curvature of $\SS_x^-$ is 
constant and equals:
$$-\frac{d}{n-1}\coth(a)+\frac{n-1-d}{n-1}\tanh(a)$$
for some $d\in\{1,\dots,n-1\}$. Finally, one observes that this quantity is
maximal when $d$ is minimal (i.e. when $d=1$). The theorem follows. 
\end{proof}

\begin{rema}
\label{rk.montiel}
The past barriers appearing in the proof are the CMC hypersurfaces presented in
Example 2 of \cite{montiel1}.
\end{rema}

By reversing the time one obtains the following result.  

\begin{teo}
\label{t.reverse-barriers-dS}
Let ${\mathcal B}^-_0(S)$ be a past complete $\mbox{dS}_n$ regular
domain, and $\widehat{\tau}: {\mathcal B}^-_0(S) \rightarrow (0, +\infty)$ 
be the reverse cosmological time function associated to ${\mathcal
  B}^-_0(S)$. 
  Then, for every $a\in (0,+\infty)$, 
the generalized mean curvature of the level set $\widehat{S_a} =
\widehat{\tau}^{-1}(a)$ admits the following estimates:
$$
\frac{1}{n-1}\coth(a) + \frac{n-2}{n-1}\tanh(a)\leq H_{\widehat{S_a}}
\leq \coth(a). 
$$
\fin
\end{teo}

\section{CMC time functions in de Sitter spacetimes}
\label{proof.dS}

In this section we prove Theorems~\ref{t.para-dS-case}
and \ref{t.main-dS-case}, and discuss CMC foliations
in elliptic de Sitter spacetimes.
The existence problem of CMC-times or CMC-foliations splits
into several  cases (essentially three) and subcases. 

\subsection{The hyperbolic case}
\label{sub.dscmc}

The proof of Theorem \ref{t.main-dS-case} is very similar 
to that  of  Theorem \ref{t.main-flat-case}. 
The only difference is that, in the de Sitter case, 
the cosmological time function does not provide a sequence 
of future asymptotic barriers (except in dimension   $2+1$).

\begin{proof}[Proof of Theorem \ref{t.main-dS-case}]
Let $(M,g)$ be a past 
incomplete $n$-dimen\-sional MGHC
spacetime  of  the de Sitter type. According to Theorem
\ref{teo.dscompact}, $(M,g)$ is the quotient of a regular domain
${\mathcal B}^+_0(S)$ by a torsion-free discrete group
$\Gamma\subset\mbox{Isom}(\mbox{dS}_n)$.   
The cosmological time $\tau:{\mathcal B}^+_0(S)\rightarrow
(0,+\infty)$ is well-defined and regular.

For every $a \in [0,+\infty]$, let $S_a=\tau^{-1}(a)$ 
and $\Sigma_a$ be the projection  of $S_a$ in 
$M\equiv\Gamma\setminus {\mathcal B}^+_0(S)$. As every compact level
set of a time function, $\Sigma_a$ is a topological 
Cauchy hypersurface in $M$ for every $a$. 
Theorem~\ref{t.barriers-dS} implies that, for every $a\in
(0,+\infty)$, the generalized mean curvature of $\Sigma_a$ satisfies
$$
-\coth(a)\leq H_{\Sigma_a}\leq -\frac{1}{n-1}\coth(a) -\frac{n-2}{n-1}\tanh(a).
$$
  
Let $(a_m)_{m\in\NN}$ be a decreasing sequence of 
positive real numbers such that $a_m\to 0$ when 
$m\to +\infty$. Observe that 
$$
-\frac{1}{n-1}\coth(a_m)-\frac{n-2}{n-1}\tanh(a_m)\to -\infty \quad\mbox{ when
}m\to\infty.
$$ 
Hence $(\Sigma_{a_m})_{m\in\NN}$ is a sequence of past asymptotic
$\alpha$-barrier in $M$ for $\alpha=-\infty$. Hence
Theorem \ref{t.half-theorem} implies that $M$ admits a partially
defined CMC-time $\tau_{cmc}:U\to (-\infty,\beta)$ where $U$
is a neighbourhood of the past end of $M$.

\subsubsection{The three-dimensional case} 
Assume $n=3$. Consider a sequence  $(b_m)_{m\in\NN}$ of increasing
positive real numbers such that $b_m\to +\infty$ when $m\to
+\infty$. For every $m\in\NN$, one has 
$$
-\coth(b_m)<-\frac{1}{2}\coth(b_m)-\frac{1}{2}\tanh(b_m)<-1
$$ 
and 
$$
-\coth(b_m) \to-1\quad\mbox{ when }m\to\infty.
$$
Hence $(\Sigma_{b_m})_{m\in\NN}$ is a sequence of future asymptotic
$\beta$-barrier in $M$ for $\beta=-1$. Therefore, Theorem \ref{t.foliation}
implies that $M$ admits admits a {\it globally defined} CMC time function
$\tau_{ cmc}:M\to (-\infty,-1)$. 

\begin{rema}
\label{rk.dsnocmc}
This argument fails if $n>3$. The problem is that the quantity 
$$-\frac{1}{n-1}\coth(a)-\frac{n-2}{n-1}\tanh(a)$$
becomes bigger than $-1$ when $a$ is large. See
\S\ref{sss.no-CMC-time} below.  
\end{rema}

\subsubsection{The almost-fuchsian case} 
In the almost-fuchsian case there is an embedded Cauchy surface
$\Sigma$ in $(M,g)$ with all principal eigenvalues  $ < -1$.
Reversing the time if needed, we can assume that $M$ is future complete.
Denote by $\Sigma_t$ the image of the hypersurface $\Sigma$ under the
time $t$ map of the Gauss flow, i.e. obtained by pushing
$\Sigma$ during a time $t$ along its normal geodesics. It is easy to
describe in our context these hypersurfaces: 
let $\widetilde{\Sigma}$ be the universal covering of $\Sigma$: 
the embedding $\Sigma \subset M$ lifts to an embedding 
$u: \widetilde{\Sigma} \rightarrow {\mathcal B}_0^+(S)$. For every $x$
in $\widetilde{\Sigma}$, there exists a unique element $u^*(x)$ of
${\mathbb H}^n_+$ such that the line $\RR.u^*(x)$ is the
$Q_{1,n}$-orthogonal of $\RR.H(x)$ where $H(x)$ is the tangent
projective hyperplane of $\widetilde{\Sigma}$ at $x$. 
In other words, we have two maps $u,u^\ast: \widetilde{\Sigma}
\rightarrow {\mathbb R}^{1,n}$ such that, for 
every $x$ in $\widetilde{\Sigma}$,
\begin{itemize}
\item $Q_{1,n}(u(x)) = 1$,
\item $Q_{1,n}(u^\ast(x)) = -1$,
\item $\langle u(x) \mid u^\ast(x) \rangle = 0$,
\item for every tangent vector $\partial_x$ at $u(x)$ we have $\langle
   u^\ast(x) \mid \partial_xu \rangle = 0$. 
\end{itemize}
Then for every $x$ in $\widetilde{\Sigma}$ we have 
$\langle u(x) \mid \partial_xu^\ast \rangle = 0$. The Weingarten operator 
for $\widetilde{\Sigma}$ is the linear operator $B$ such that
$B(\partial_xu) = -\partial_xu^\ast$ for every tangent vector $\partial_x$.

The Gauss flow is described as follows: for every $t\geq0$, let 
$u_t: \widetilde{\Sigma} \rightarrow \mbox{dS}_n \subset {\mathbb
  R}^{1,n}$ defined 
by $u_t(x) = \cosh(t)u(x) + \sinh(t)u^\ast(x)$. 
Observe that since we have selected $u^\ast(x)$ in ${\mathbb H}_+^n$ the 
$u_t(x)$ (for a fixed $x$) describes a future oriented geodesic ray
starting from $u(x)$. The projection in $M$ of the image
$\widetilde{\Sigma}_t$ of $u_t$ is the hypersurface $\Sigma_t$. 

For a fixed $t$, the differential of
$u_t$ evaluated on a tangent vector $\partial_x$ is 
$\cosh(t)\partial_xu + \sinh(t)\partial_xu^\ast =
(\cosh(t)Id - \sinh(t)B)(\partial_xu)$.

By assumption,  the principal curvatures of $\Sigma$, i.e. the
eigenvalues of $B$, 
are less than $-1$. It follows that $u_t$ is an immersion for every 
$t\geq0$: the Gauss flow is defined for all positive $t$.
Moreover, the differential of $u^\ast_t$ evaluated on $\partial_x$
is $(\sinh(t)Id - \cosh(t)B)(\partial_xu)$. It follows that the
Weingarten operator 
for $B_t$ is $-(\tanh(t)Id - B)(Id - \tanh(t)B)^{-1}$. In particular, the mean
curvature of $\Sigma_t$ is smaller than $-1$ for every $t\geq 0$, and
tends to $-1$ when $t\rightarrow +\infty$. 

Now, we claim that given an increasing
sequence $(t_m)_{m\in\NN}$ of real numbers such that $t_m\to\infty$
when $m\to \infty$, the sequence of hypersurfaces
$(\Sigma_{t_m})_{m\in\NN}$ is a sequence of future asymptotic
$\beta$-barrier in $M$ for $\beta=-1$.  The only remaining point
to check is that $(\Sigma_{t_m})_{m\in\NN}$ tends to the future end
of $M$ when $m\to+\infty$. But this is clear: let $T_0$ be the minimal value
of the cosmological time function on $\Sigma$. Then the cosmological 
time function restricted to $\Sigma_t$ is everywhere bigger than
$T_0+t$. The claim follows.

Hence Theorem \ref{t.foliation} implies that $M$ admits a globally
defined CMC-time $\tau_{cmc}:M\to (-\infty,-1)$. 
\end{proof}

\begin{rema}
\label{rk.fuchsianalmost}

We define (future complete) \emph{fuchsian} de Sitter spacetimes as MGHC
de Sitter spacetimes $(M,g_0) = \cB_0^+(S)$ where the M\"{o}bius manifold $S$

is a quotient $\Gamma\backslash{U}$ of a proper round ball $U$ in $\SS_+^{n-1}$. 
The metric of a Fuchsian spacetime is a warped product of the form  
 $-dt^2  + w(t)^2 h$, where $h$ is independent of $t$. 
Any metric of this form 
admits a timelike homothety and is conformal to a static spacetime.

Observe that in particular 
the holonomy group $\Gamma$ is conjugate in $\mbox{SO}_0(1,n)$
to a lattice of $\mbox{SO}_0(1,n-1)$;
$\Gamma$ preserves a totally geodesic
hypersurface $\HH^{n-1}$ in $\HH^n$.

We claim that $(M, g_0)$ is almost-fuchsian.
To see this, consider a hypersurface $\Sigma$ dual to a hypersurface
in $\HH^n$ all the principal curvatures of which are very small (this
last hypersurface can be obtained by taking the image of the totally
geodesic hypersurface $\HH^{n-1}$ under the time $t$ map of the Gauss flow for $t$
small). 

If another Lorentz metric  $g$  of dS type is a small
deformation of the fuchsian metric $g_0$, then the hypersurface $\Sigma$
also has all its principal curvatures $ < -1$ (with respect to  $g$). 
\end{rema}

\begin{rema}
In  dimension $2+1$,  Theorem \ref{t.main-dS-case} can also be
deduced from the existence of foliation of hyperbolic ends by surfaces
with constant Gauss curvature (see~\cite{Bar.Zeg}).  
\end{rema}

\subsubsection{A regular spacetime with no CMC time function}
\label{sss.no-CMC-time}
For every $n\geq 4$, there exists $n$-dimensional MGHC regualr
spacetimes that do not admit any CMC time function. Here a
construction of such a spacetime. Let $n\geq 4$ and choose as
M\"{o}bius surface $S$ the complement in ${\mathbb S}_+^{n-1}$ of two 
points, say $p_1$ and $p_2$. Let $P_1$ and $P_2$ be the 
projective hyperplanes in $\SS(\RR^{n+1})$ which are tangent to
$\SS_+^{n-1}$ respectively at $p_1$ and $p_2$. The intersection 
$Q=P_1\cap P_2$ is a spacelike totally geodesic subspace
of dimension $n-2$ in $\DS_n$, homeomorphic to $\SS^{n-2}$. The domain
$\BB_0^+(S)$ is by 
definition the intersection of the futures of $P_1$ and the future of
$P_2$. It can be easily proved that the cosmological time
function $\wt\tau$ of ${\mathcal B}_0^+(S)$ is just the lorentzian
distance to the spacelike totally geodesic $(n-2)$-sphere $Q$. 
Using this, one can verify that, for every $a$, the level set
$S_a=\wt\tau^{-1}(a)$ is a Cauchy hypersurface in $\BB_0^+(S)$ which is
homeomorphic to the $\SS^{n-2}\times \RR$, and has constant mean
curvature equal to 
$$
-\frac{1}{n-1}\coth(a)-\frac{n-2}{n-1}\tanh(a)
$$
(the calculation of the mean curvature is entirely similar to the
estimates of the curvature of the hypersurface $\SS^+$ in the proofs
of Theorem~\ref{t.barriers-AdS} and~\ref{t.barriers-dS}).
Now, observe that the regular domain ${\mathcal B}_0^+(S)$
admits (regular) Cauchy compact quotients: if $\Gamma$ is a cyclic
group generated by a hyperbolic element of $\mbox{SO}_0(1,n)$ fixing
the points $p_1$ and $p_2$, then $\Gamma$ acts properly
discontinuously on $\BB_0^+(S)$ and the projection $\Sigma_a$ of $S_a$
in the quotient $M:=\Gamma\setminus B_0^+(S)$ is a Cauchy hypersurface
homeomorphic to $\SS^{n-2}\times\SS^1$. Moreover, for every $a$, the
hypersurface $\Sigma$ has  constant mean curvature
equal to $-\frac{1}{n-1}\coth(a)-\frac{n-2}{n-1}\tanh(a)$. 
Hence $\cF=\{\Sigma_a\}_{a\in (0,+\infty)}$ is a CMC foliation of
$M$. But the mean curvature of the leafs of $\cF$ is not
monotonous (it increases for $a$ small, but decreases for $a$ large).
In particular, $M$ does not admit any CMC time function (if there
would exist a CMC time function, then the hypersurface $\Sigma_a$
would be a fiber of this CMC time function for every $a$, and thus,
the mean curvature of $\Sigma_a$ would be a monotonous function of $a$).

This raises the following question.

\begin{question}
Do every MGHC regular spacetime admit a global CMC foliation with
compact leaves?
\end{question}

\subsection{The elliptic case}
\label{sub.ellicmc}
\subsubsection{de Sitter space}
We first consider the case of de Sitter space itself $\mbox{dS}_n$.
A key fact is that compact CMC hypersurfaces in $\mbox{dS}_{n}$ are umbilical
(see \cite{montiel2}; this is of course reminiscent of Alexandrov
rigidity theorem which states that is any compact CMC hypersurface in
the Euclidean space is a round sphere). More precisely, they are the
intersections between $\mbox{dS}_{n} = \{ Q_{1,n} = 1 \}$ and the
affine spacelike hyperplanes of the Minkowski space ${\mathbb
  R}^{1,n}$. Such an hyperplane is defined as the set $H_{(t, v)} = \{
x / \langle x \mid v \rangle = \sinh(t) \}$ where $v$ is a vector of
norm $-1$ in the future cone of the Minkowski space, i.e. an element
of the hyperbolic space ${\mathbb H}^{n} = \{ Q_{1,n} = -1 \}$, 
and $t$ a real number. Then, the intersection $S_{(t, v)} =
H_{(t,v_{0})} \cap \mbox{dS}_{n}$ is an umbilical sphere, and every
closed CMC surface in $\mbox{dS}_{n}$ must be such an intersection. 
In other words, ${\mathbb H}^{n} \times {\mathbb R}$ is the space of
umbilical spheres.

The mean curvature of $S_{(t,v)}$ is $-\tanh(t)$. It follows that if
$S_{(t, v)}$ is in the future of $S_{(t',v')}$,  then the mean
curvature of the former is less than the mean curvature of the later. 
This phenomenom is actually valid locally.

\begin{lem} 
 Let $U$ be an open subset of $\mbox{dS}_n$ endowed with an umbilical
 foliation ${\mathcal F}$  with compact leaves.  Then, the mean
 curvature function of ${\mathcal F}$ is decreasing. In particular,
 $\mbox{dS}_n$ has no CMC time. 
\end{lem}

\begin{proof}
By contradiction, assume that the mean curvature is somewhere
increasing (or just non-decreasing). This will be true on an open
${\mathcal F}$-saturated set, we can thus assume that this  holds  on
all $U$.  Therefore, on $U$, we have a CMC time. By a well known
property, any other compact CMC hypersurface in $U$ is a leaf of
${\mathcal F}$. This is obviously false:   take $S$ a leaf of
${\mathcal F}$, and $S^\prime$ an umbilical hypersurface close to it,
then $S^\prime$ will be contained in $U$, but is not necessarily a
leaf of $\mathcal F$. 

 Observe in fact that for a global foliation of $\mbox{dS}_n$, leaves
 accumulate to the two  boundary components, which can be thus  seen
 as umbilical hypersurfaces, but with infinite curvature. More
 formally, the curvature of leaves decreases (with time) from $+
 \infty$ to $-\infty$. 
 \end{proof}

We want to describe now CMC-foliations in $\mbox{dS}_{n}$. The
following Proposition gives a complete description.

\begin{prop}
\label{p.foliations-de Sitter}
There is a 1-1 correspondance between CMC-foliations with compact
leaves in $\mbox{dS}_{n}$ and inextendible timelike curves in
${\mathbb H}^{n} \times {\mathbb R}$ equipped 
with the lorentzian metric $ds_{hyp}^{2} - dt^{2}$ where
$ds_{hyp}^{2}$ is the hyperbolic metric of ${\mathbb H}^{n}$.
\end{prop}

\begin{proof}
Let $\mathcal F$ be a CMC-foliation with compact leaves. In order to
simplify the proof, we assume 
that $\mathcal F$ is $C^{1}$, but see remark~\ref{rk.feuillise}.
The leaves are umbilical spheres $S_{(t,v)}$. Observe that since the
leaves are disjoint one to the other, two different leaves must have different
parameter $t$. By Reeb stability theorem (see \cite{godbillon}), since
every leaf is a sphere, the foliation is trivial: there is a map 
$f: \mbox{dS}_{n} \to {\mathbb R}$ such that the leaves of $\mathcal F$ 
are the fibers of $f$. It follows that there is a curve 
$c_{\mathcal F}: I \rightarrow {\mathbb H}^{n} \times {\mathbb R}$
such that the leaves of $\mathcal F$ are the umbilical spheres
$S_{(t(s), v(s))}$ where $I \subset \mathbb R$  
and $c_{\mathcal F}(s) = (t(s), v(s))$. Since the map $s \to t(s)$ is
1-1, we can choose that the parameter $s$ 
so that $t(s) = s$, i.e. we can parametrize $c_{\mathcal F}$ by the
first factor $t$. 

Consider any $C^{1} $curve $c: I \rightarrow {\mathbb H}^{n} \times
{\mathbb R}$: the umbilical spheres $S_{c(t)}$ may be non-disjoint. We make
the following

\begin{claim} 
The spheres $S_{c(t)}$ are pairwise disjoint if and only
  if tangent vectors $v'(t)$ have hyperbolic norm less than $1$.
\end{claim}

We first consider the case $n=1$. Then $v(t) = (\sinh(\eta(t)),
\cosh(\eta(t)))$ where $t \to \eta(t)$ is a $C^{1}$ map. The elements
of the $0$-sphere $S_{(t,v(t))}$ are 
$(\cosh(a), \sinh(a))$ and $(-\cosh(b), \sinh(b))$ where $a$, $b$ satisfy:

\begin{eqnarray*}
\cosh(a)\sinh(\eta) - \sinh(a)\cosh(\eta) = \sinh(t)\\
-\cosh(b)\sinh(\eta) - \sinh(b)\cosh(\eta) = \sinh(t)
\end{eqnarray*}

Hence, we have $a = t - \eta$ and $b = t + \eta$. But the $0$-spheres
$S_{c(t)}$ are disjoint if and only if the maps $t \to a$ and $t \to
b$ are increasing. This is equivalent to the absolute value of
$\eta'(t)$ being strictly less than $1$. The claim follows since the
hyperbolic metric of ${\mathbb H}^{1}$ is $d\eta^{2}$.

Assume now $n \geq 2$. Let $P$ be any $2$-plane in ${\mathbb R}^{1,n}$
on which the restriction of $Q_{1,n}$ has signature $(1,1)$. Let
$\pi_{P}: {\mathbb R}^{1,n} \to P$ be the orthogonal projection. 
If the $S_{c(t)}$ are two by two disjoint the same is true
for the intersections $P \cap S_{c(t)}$, and conversely, if $P \cap
S_{c(t)}$ and $P \cap S_{c(t')}$ are disjoint for every $2$-plane as
above, then $S_{c(t)}$ and $S_{c(t')}$ are disjoint. Now observe that
the intersection $P \cap S_{c(t)}$ is nothing but the set of points 
$x$ in $P \cap \mbox{dS}_{n} \approx \mbox{dS}_{1}$ satisfying
$\langle x \mid \pi_{P}(v) \rangle = \sinh(t)$. 
Hence, since the $n=1$ case has been proved, the spheres $S_{c(t)}$
are all disjoint if and only if for every $2$-plane $P$ as above the
norm of $d\pi_{P}(v'(t))$ is less than one. But, using the natural
parallelism of ${\mathbb R}^{1,n}$, the spacelike vector $v'(t)$ has
Minkowski norm less than $1$ if and only if all the vectors
$d\pi_{P}(v'(t)) = \pi_{P}(v'(t))$ have Minkowski norm less than
$1$. The claim follows. 

According to the claim, the curve $c_{\mathcal F}: I \to \mathbb R$ is
a timelike curve in ${\mathbb H}^{n} \times \mathbb R$. If this curve
is extendible, then it means that some umbilical curve $S_{(T,V)}$ is
disjoint from all the $S_{c_{\mathcal F}(t)}$. This is a contradiction
since $\mathcal F$ foliates the entire de Sitter space.  Hence,
$c_{\mathcal F}$ is inextendible. 

Conversely, for every inextendible timelike curve $c$ in $\HH^{n}
  \times \RR$, the arguments above show that $t\to S_{c(t)}$ is a
$1$-parameter family of umbilical spheres which are pairwise
  disjoint. Since the projection on the second  factor of is a Cauchy
  time function on the globally hyperbolic space ${\mathbb H}^{n}
  \times \mathbb R$, the mean curvature $t$ must takes all value in
  $]-\infty, +\infty[$. We leave to the reader the proof that the
  continuity of $c$ implies that the spheres $S_{c(t)}$ cover all the 
de Sitter space. It follows that the spheres that they are the leaves of a
  CMC-foliation $\cF_c$. 
\end{proof}

\begin{cor}
\label{cor.CMCds}
There are infinitely many non-isometric CMC-foliations of the
de Sitter space $\dS_n$. \fin  
\end{cor}

\begin{remark}
\label{rk.feuillise}
\begin{enumerate}
\item Proposition~\ref{p.foliations-de Sitter} actually shows that the
  modulus space of CMC foliations of the de Sitter space $\dS_n$ up to
  isometry is enormous: this is an open set in an infinite
  dimensional vector space. 
\item 
Proposition~\ref{p.foliations-de Sitter} provides many examples of CMC
foliations of $\dS_n$ with poor regularity. Indeed, consider a inextendible
timelike curve $c$ in ${\mathbb H}^{n} \times {\mathbb R}$ (equipped with the
lorentzian metric $ds_{hyp}^{2} - dt^{2}$). The proof of
Proposition~\ref{p.foliations-de Sitter} shows how to associate with the curve
$c$ a CMC foliation $\cF_c$ of $\dS_n$. Each leaf of the foliation $\cF_c$ is
an  umbilical sphere in $\dS_n$; in particular, it is an analytic submanifold
of $\dS_n$. Nevertheless, it follows easily from the construction that the
tranversal regularity of the foliation $\cF_c$ is exactly the same as the
regularity of the curve $c$. More precisely, if $\gamma$ is analytic curve
tranversal to the foliation $\cF_c$, the tangent plane of the leaves of
$\cF_c$ varies in a $C^k$ way along $\gamma$ if and only if the curve $c$ is
$C^k$. Therefore, a curve $c$ which is $C^k$ but not $C^{k+1}$ yields a CMC
foliation $\cF_c$ of $\dS_n$ which is $C^k$ but not $C^{k+1}$.
\item It is well-known that the notion of timelike curve in a lorentz
manifold extend to the non-differentiable case: here, it can be defined as
curves $c: t \to {\mathbb H}^{n} \times \mathbb R$ such that $c(t)$ is
in the strict future of $c(t')$ for all real numbers $t' < t$. Such
curves are automatically Lipschitz (see \cite{beem}). It is quite obvious that
timelike curves in this more general meaning also provide
CMC-foliations which are only Lipschitz regular. 
\item  In Proposition~\ref{p.foliations-de Sitter}, we only considered
  foliations with compact leaves. It is suggestive
to relax this condition, i.e. to ask whether CMC-foliations  with non
compact leaves of $\mbox{dS}_n$  exist and how they behave?
\item The opposite of the mean curvature  of an umbilical foliation is a 
time function. But, not all umbilical time functions are equally 
``tame''.  For instance, given any (spacelike compact) hypersurface 
$S$ in $\mbox{dS}_n$, its isometry group $G_S$  (i.e. isometries of
$\mbox{dS}_n$ preserving it) has umbilical orbits.  The so-obtained
time is $G_S$-invariant.  No other  time function can have a
``comparable'' symmetry group. It is interesting to characterize,
variationally, say, these extra-symmetric time functions.
\end{enumerate}
\end{remark}

\subsubsection{Non-trivial quotients of $\mbox{dS}_{n}$}
In general, an elliptic MGHC de Sitter spacetime is the quotient
of $\mbox{dS}_n$ by a \emph{finite\/} group $\Gamma$ acting freely on
$\mbox{dS}_{n}$. The group $\Gamma$ admits a fixed point $v_{0}$ in
${\mathbb H}^{n}$. For every real number $t$, the umbilical sphere
$S_{(t,v_{0})}$ is preserved by $\Gamma$: it projects in the quotient
$M = \Gamma\backslash\mbox{dS}_{n}$ on a umbilical
hypersurface. Hence, varying $t$, we obtain a CMC foliation ${\mathcal
  F}_{0}$ in $M$. Observe that $M$ admits no CMC time function, since such a
CMC time function would lift in $\mbox{dS}_n$ to a CMC time
function. Furthermore: 

\begin{lem}
\label{le.CMCds}
Every compact CMC hypersurface in $M$ is a leaf of ${\mathcal F}_{0}$.
\end{lem}

\begin{proof}
Let $S$ be a CMC hypersurface in $M$. It lifts to a compact CMC
hypersurface in $\mbox{dS}_{n}$, i.e. to some umbilical sphere
$S_{(t,v)}$. It is easy to show that for any isometry $\gamma$ of 
$\mbox{dS}_{n}$, either we have $\gamma S_{(t,v)} = S_{(t,v)}$, or
there is a transverse intersection between $\gamma S_{(t,v)}$ and
$S_{(t,v)}$. Since here $S_{(t,v)}$ is the lifting of $S$, the former
case cannot occur when $\gamma$ belongs to $\Gamma$. Hence, $v$ must
be a fixed point of $\Gamma$. Assume $v \neq v_{0}$. Then,
$S_{(0,v_{0})}$ is the unit sphere in the euclidean space
$v_{0}^{\perp} \approx {\mathbb R}^{n}$, and $v_{1}^{\perp} \cap
v_{0}^{\perp}$ is a  $\Gamma$-hyperplane in this euclidean space. The
orthogonal to this hyperplane for the euclidean metric in
$v_{0}^{\perp}$ intersects the unit sphere in two points which are
both fixed by $\Gamma$ (indeed, these points are fixed individually
and not permuted, since one of them belongs to the future of
$v_{1}^{\perp}$ in ${\mathbb R}^{{1,n}}$ and the other belongs to the
past of $v_{1}^{\perp}$). This is a contradiction since the action of
$\Gamma$ on $\mbox{dS}_{n}$ is free. Hence, $v = v_{0}$: the
hypersurface $S$ is a leaf of ${\mathcal F}_{0}$. 
\end{proof}

Corollary~\ref{cor.CMCds} and Lemma~\ref{le.CMCds} give the proof
of Theorem~\ref{t.elli-dS-case}.

\subsection{The parabolic case}
\label{sub.paracmc}
Consider a parabolic standard spacetime ${\mathcal B}_0^+(S)$. By
 definition of parabolic spacetimes, $S$ is the sphere ${\mathbb
 S}^{n-1}_+$ of one point $r_0$. The hyperbolic space ${\mathbb
 H}^{n}_+$ is foliated by umbilical hypersurfaces with constant mean
 curvatures $-1$: the horospheres based at $r_0$. The dual to these
 hypersurfaces are umbilical hypersurfaces with the same constant
 mean curvature $-1$, and  foliate ${\mathcal B}_0^+(S)$ (these
 hypersurfaces are not umbilical spheres, but it is not a
 contradiction with Montiel's theorem since they are not compact!).  
It follows that ${\mathcal B}_0^+(S)$  admits no CMC time function 
(since as explained above, if such a CMC time function would exist,
 then any CMC hypersurface would be a level set of this function; in
 particular, there would exist at most one CMC hypersurface with mean
 curvature $-1$ in ${\mathcal B}_0^+(S)$). 

Every future complete parabolic  MGHC dS spacetime is a quotient 
$M = \Gamma\setminus {\mathcal B}_0^+(S)$ where $\Gamma$ is a subgroup
of $\mbox{SO}_0(1,n)$ preserving $\infty$.  As in previous case, we
have a CMC-foliation but no CMC-time. Moreover, let $\Sigma$ be any
closed CMC hypersurface. It is tangent to two leaves of the
CMC-foliation, one of these leaves being in the future of $\Sigma$,
and the other in the past. By the maximum principle, $\Sigma$ has mean
curvature $-1$; by the equality case of the maximum principle it
follows that $\Sigma$ is equal to the CMC-leaves. In particular, the
CMC-foliation is unique. This completes the proof of Proposition~\ref{t.para-dS-case}. $\Box$

\begin{rema}
Proposition~\ref{t.para-dS-case} also follows directly from \cite{montiel2}.
\end{rema}

\end{document}